\RequirePackage{fix-cm}
\documentclass[envcountsect]{svjour3}
\smartqed  % flush right qed marks, e.g. at end of proof
\usepackage{amsmath,amssymb,color}
\usepackage[english]{babel}
\usepackage[utf8]{inputenc}
\usepackage[T1]{fontenc}
\usepackage{authblk,mathrsfs}
\usepackage{algorithm,algpseudocode}
\usepackage{graphicx,subcaption}
\usepackage{hyperref}
\usepackage{booktabs}
\usepackage{caption}
\usepackage{subcaption}
\usepackage{pdfsync}
\newcommand{\RR}{\mathbb{R}}
\newcommand{\CC}{\mathbb{C}}
\newcommand{\mrm}{\mathrm}
\newcommand{\mL}{\mathrm{L}}
\newcommand{\mH}{\mathrm{H}}
\newcommand{\mJ}{\mathrm{J}}
\newcommand{\mV}{\mathrm{V}}
\newcommand{\Id}{\mathrm{Id}}

\newcommand{\mS}{\mathrm{S}}

\newcommand{\mX}{\mathrm{X}}
\newcommand{\mbH}{\mathbb{H}}
\newcommand{\mbV}{\mathbb{V}}

\newcommand{\Vh}{\mathrm{V}_{h}}

\newcommand{\Th}{t_h}
\newcommand{\rh}{\rho_h}

\newcommand{\mbVh}{\mathbb{V}_{h}}

\newcommand{\rref}{\textsc{r}}

\newcommand{\by}{\boldsymbol{y}}
\newcommand{\bx}{\boldsymbol{x}}

\newcommand{\bn}{\boldsymbol{n}}

\newcommand{\bfp}{\mathbf{p}}
\newcommand{\bfff}{\mathbf{f}}

\newcommand{\bfu}{\mathbf{u}}
\newcommand{\bfv}{\mathbf{v}}
\newcommand{\bfq}{\mathbf{q}}
\newcommand{\bfg}{\mathbf{g}}
\newcommand{\bfb}{\mathbf{b}}
\newcommand{\dir}{\textsc{d}}
\newcommand{\neu}{\textsc{n}}

\newcommand{\lbr}{\lbrack}
\newcommand{\rbr}{\rbrack}
\newcommand{\ctru}{\mathfrak{u}}
\newcommand{\ctrv}{\mathfrak{v}}
\newcommand{\ctrw}{\mathfrak{w}}
\newcommand{\ctrp}{\mathfrak{p}}
\newcommand{\ctrq}{\mathfrak{q}}
\newcommand{\ctrf}{\mathfrak{f}}
\renewcommand{\dim}{\mrm{dim}}
\newcommand{\bfA}{\mathbf{A}}
\newcommand{\bfB}{\mathbf{B}}
\newcommand{\bfT}{\mathbf{T}}

\newcommand{\bfQ}{\mathbf{Q}}
\newcommand{\iu}{i}

\spnewtheorem{thm}{Theorem}[section]{\bfseries}{\itshape}
\spnewtheorem{cor}[thm]{Corollary}{\bfseries}{\itshape}
\spnewtheorem{lem}[thm]{Lemma}{\bfseries}{\itshape}
\spnewtheorem{prop}[thm]{Proposition}{\bfseries}{\itshape}
\spnewtheorem{defn}[thm]{Definition}{\bfseries}{\itshape}% maybe \upshape?
\spnewtheorem{as}[thm]{Assumption}{\bfseries}{\itshape}
\spnewtheorem{rmrk}[thm]{Remark}{\bfseries}{\itshape}
\spnewtheorem{exmpl}[thm]{Example}{\bfseries}{\itshape}
\spnewtheorem{cond}[thm]{Condition}{\bfseries}{\itshape}

\journalname{Numerische Mathematik}
\begin{document}

\title{Robust treatment of cross-points in Optimized Schwarz Methods%\thanks{Grants or other notes
}


\author{Xavier Claeys         \and
        Emile Parolin %etc.
}


\institute{X. Claeys \at
              Sorbonne Universit\'e, Université Paris-Diderot SPC, CNRS, Inria,
  Laboratoire Jacques-Louis Lions, \'equipe Alpines\\
              \email{claeys@ljll.math.upmc.fr}           %  \\
           \and
           E. Parolin \at
              POems, UMR CNRS/ENSTA/Inria, Institut Polytechnique de Paris\\
              \email{emile.parolin@inria.fr}
}

\date{Received: date / Accepted: date}
\date{}

\maketitle

\begin{abstract}
    In the field of Domain Decomposition (DD), Optimized Schwarz Method (OSM)
    appears to be one of the prominent techniques to solve large scale
    time-harmonic wave propagation problems. It is based on appropriate
    transmission conditions using carefully designed impedance operators to
    exchange information between subdomains. The efficiency of such methods
    is however hindered by the presence of cross-points, where more than two
    subdomains abut, if no appropriate treatment is provided.

    In this work, we propose a new treatment of the cross-point issue for the
    Helmholtz equation that remains valid in any geometrical
    interface configuration. We exploit the multi-trace formalism to define a
    new exchange operator with suitable continuity and isometry properties.
    We then develop a complete theoretical framework that generalizes
    classical OSM to partitions with cross-points and contains a rigorous proof
    of geometric convergence, uniform with respect to the mesh discretization, for
    appropriate positive impedance operators. Extensive numerical results in 2D and
    3D are provided as an illustration of the proposed method. 
%Insert your abstract here. Include keywords, PACS and mathematical
%subject classification numbers as needed.
\end{abstract}

%\section{Introduction}
%\label{intro}
%Your text comes here. Separate text sections with
%\section{Section title}
%\label{sec:1}
%Text with citations \cite{RefB} and \cite{RefJ}.
%\subsection{Subsection title}
%\label{sec:2}
%as required. Don't forget to give each section
%and subsection a unique label (see Sect.~\ref{sec:1}).
%\paragraph{Paragraph headings} Use paragraph headings as needed.
%\begin{equation}
%\end{equation}
%
%% For one-column wide figures use
%\begin{figure}
%% Use the relevant command to insert your figure file.
%% For example, with the graphicx package use
%  \includegraphics{example.eps}
%\end{figure}
%%
%% For two-column wide figures use
%\begin{figure*}
%% Use the relevant command to insert your figure file.
%% For example, with the graphicx package use
%  \includegraphics[width=0.75\textwidth]{example.eps}
%\end{figure*}
%\begin{table}
%\begin{tabular}{lll}
%\hline\noalign{\smallskip}
%\end{tabular}
%\end{table}

\section*{Introduction}

Domain Decomposition (DD) for time-harmonic wave propagation is presently an
active field of research as the numerical simulation of large scale problems
remains a challenge in scientific computing. The additional difficulty of such
problems, in comparison to elliptic problems, mainly lies in the (a priori) indefiniteness
of the Helmholtz equation and related linear systems after discretization.
In this work, we are interested in the sub-class of so-called (non-overlapping)
Optimized Schwarz Methods (OSM) that appear as the most established approach
in a wave context. In the context of waves, OSM dates back to the PhD thesis of
Despr\'{e}s~\cite{MR1071633,MR1105979,MR1291197,MR1227838} where the idea to
use Robin or impedance like transmission quantities was first introduced and a
proof of algebraic convergence using energy estimates was derived.  Many
refinements over this initial idea were proposed since then in order to improve
the rate of convergence exhibited.
Most of the methods that were latter derived rely on
the definition of a generalized Robin quantity with the introduction of an
impedance operator.  In this spirit, second order impedance operators (which
possess sufficient properties for guaranteed convergence) were
introduced by Gander, Magoul\`es and Nataf~\cite{MR1924414} and detailed in later
works~\cite{MR3549886,MR2097757}. A popular approach was then developed which
consists in using a (high-order) absorbing boundary condition (ABC) as
transmission condition. The underlying idea is to approximate the Dirichlet-to-Neumann
(DtN) map in the complementary domain, which provides exact transparent
conditions albeit at a prohibitive numerical cost.
This is the approach adopted in~\cite{Boubendir2012,ElBouajaji2015} for
instance where Pad\'e-approximants of the square root operator are used to construct
high-order ABC.\@ However, despite their efficiency in practice, such methods
lack a rigorous analysis of convergence since only partial proofs in specific
geometries are available.
Another alternative, advocated
in~\cite{MR1764190,LecouvezCollinoJoly,lecouvez:tel-01444540}, is to use
suitable non-local operators, realized in practice with integral operators, as
impedance operators.
One of the strengths of this later approach is to rely on a solid theoretical basis
that systematically guarantees (\(h\)-uniform~\cite{Claeys2019}) geometric
convergence, provided that certain properties of injectivity, surjectivity and
positivity (in suitable trace spaces) are satisfied by the impedance operator.

For realistic large scale applications, DD methods should be applicable to
domain partitions with cross-points, where at least three subdomains share a
common vertex. The presence of such points can be an issue both for the
analysis at the continuous level and in practice for numerical
implementations.
For DD methods used in conjunction with zeroth order transmission operators,
the convergence proof is established at the continuous level and in the case of
mixed finite element discretizations, for instance in~\cite{MR1291197}.
Interestingly, this particular choice of discretization avoids degrees of
freedom (DOF) at cross-points and therefore the issue altogether.
In contrast, the cross-point issue arises if one makes the choice of using nodal
finite element discretizations, see~\cite{Loisel2013,St-Cyr2009}.
Gander and Kwok~\cite{Gander2013} pointed that straightforward nodal discretization
of OSM can diverge and that the continuous proof (based on Lions' energy estimates)
fails to carry over to the discrete setting in general.
Dual-Primal treatment of the problem at the discrete level has been developed
which introduces additional global unknowns at cross-points effectively coupling
all subdomains~\cite{Boubendir2003,Bendali2006}. This leads to a global
indefinite system that needs to be solved at each iteration.
As regards the continuous theory available for DD methods constructed using
non-local operators~\cite{MR1764190,LecouvezCollinoJoly,lecouvez:tel-01444540},
it rests unfortunately on the strong hypothesis of the absence of cross-points
between interfaces~\cite[Rem.  3]{lecouvez:tel-01444540}.
Analysis suggests that this issue is related to the exchange operator being not
continuous in proper trace norms in the presence of cross-points.
Recently, the cross-point problem has been addressed by Despr\'es, Nicolopoulos
and Thierry~\cite{despres:hal-02612368} for a particular class of second-order
transmission operators for which convergence of iterative algorithms
is proven by using energy estimate techniques, though without any estimation of
the convergence rate. Besides, Modave \textit{et al}~\cite{Modave2020b} presented a
treatment of cross-points in the context of high-order ABC based transmission conditions.
This later approach is however only valid on Cartesian like partitions of the mesh,
allowing only cross-points where exactly four domains abut (in 2D),
and no convergence theory is provided. It is clear
that being able to deal with more general partitions, generated for instance by
graph partitioners, is a highly desirable property.

The goal of this work is to use the clean treatment of cross-points from the
Multi-Trace (MTF) formalism~\cite{MR3069956}, initially developed for
boundary integral equations (BIE), to investigate OSM.
The main idea is to introduce a regularized version of the exchange operator,
that remains isometric and continuous regardless of cross-points.
The starting point is to recognize that if one uses positive impedance
operators in the DD algorithm, it is possible to define a scalar product on the
multi-trace space (collection of traces of local solutions in the subdomains).
We then use this scalar product to define an orthogonal projector onto the
single-trace space (collection of traces that match across interfaces), which
is a closed subspace of the multi-trace space, in a very natural way.
This provides a discrete characterization of the continuity of both the
Dirichlet and Neumann traces across interfaces, that remains valid in the
presence of cross-points. The definition of the exchange operator exploits this
characterization using orthogonal projection and can be realized in practice
by solving a positive linear system posed on the skeleton of the partition.
As a result, since its computation amounts to solving a linear system, the
exchange operator is \textit{a priori} non-local.
However, the structure of the auxiliary problem and in particular
its definiteness are propitious to an efficient inversion, even in a
distributed-memory parallelization context.
Closely related ideas have been developed in a previous
work~\cite{claeys2019new} by the first author, at the continuous level.
In contrast to this work, the exchange operator was there defined explicitly
in terms of boundary integral operators.

Having defined this new way of exchanging traces across interfaces, the
classical DD algorithms remain unchanged and the coupled local systems can be
equivalently recast as a problem posed on the skeleton as usual.
In fact, in the absence of cross-point, our approach reduces to traditional OSM so
that our work can be seen as a generalization of previously described methods.
The continuity and isometric properties of the exchange operator together with
the contractivity of the local scattering operators yield immediately geometric 
convergence of the Richardson algorithm, despite the presence of cross-points.
We show that in the case of geometric partitioning independent of the 
triangulation, using uniformly bounded impedance operators with respect to mesh
discretization allows to obtain a uniform convergence rate. 

The present contribution describes what we believe is the first DDM substructuring
strategy for waves with guaranteed geometric convergence regardless of the presence
of cross-points. We also provide a theoretical framework that applies
in general geometric configurations while previous theoretical contributions
on OSM either discarded the presence of cross-points~\cite{MR1764190,LecouvezCollinoJoly,lecouvez:tel-01444540},
proposed a convergence result with no estimate of the convergence rate~\cite{MR1291197},
or considered positive definite problems~\cite{MR3519297,Loisel2013,MR3013465}.
When there is no cross-point, and for a proper choice of impedance, the DDM
approach presented here coincides with the method of Despr\'{e}s
\cite{MR1071633,MR1105979,MR1291197,MR1227838} and, as a byproduct, our analysis also
provides new convergence estimates for Despr\'{e}s algorithm, see Example
\ref{ConvergenceRateDespresImpedanceFixedPartition}.

Although derivation of high frequency robust numerical strategies is not the primary focus
of the present contribution, our analysis remains valid in the high frequency
regime. This does not mean that our method enjoys a frequency uniform convergence rate. However 
Proposition \ref{ResolventBound} establishes a link between the discrete inf-sup condition of
the initial wave propagation problem and the convergence rate of our DDM strategy and, under
the minimal assumption that the initial finite element matrix is invertible,
this correspondence holds independently of the frequency regime and the mesh width parameter.

This article is organized as follows. In Section~\ref{PbUnderStudy} we shortly describe the
problem under consideration, before describing the geometric partitioning in
Section~\ref{GeomPart} and an adapted functional framework in
Section~\ref{MultiDomainFunctionalSetting}  and \ref{sec:tracenorms}.
The definitions of impedance operators and associated
scalar products are then given in Section~\ref{sec:impedance}. In
Section~\ref{sec:reformulation}, the essential Lemma~\ref{DefOrthoProj} states
the properties of the exchange operator and the discrete characterization of
the continuity of Dirichlet and Neumann traces.
Section~\ref{sec:reduction} recasts the problem at interfaces and proves that
the equivalent problem and the local sub-problems are well posed.
Section~\ref{sec:algorithm} describes the usual Richardson algorithm of our DD method
and Theorem~\ref{ThmConvergenceEstimate} states that one obtains geometric convergence
of the iterative solution towards the solution of the original model problem.
Section~\ref{sec:discreteStability} on discrete stability gives an explicit
lower bound for the inf-sup constant associated to the problem on the skeleton.
This result yields uniform rate of convergence of iterative algorithms with
respect to the mesh discretization, in the particular case of uniformly bounded
impedance operators, as stated  in Corollary~\ref{UniformVolumicConvergenceRate},
Section~\ref{sec:fixedpartition}. We show how the proposed method is a generalization
of classical OSM in Section~\ref{sec:generalization}.
The algorithm in matrix form is detailed in Section~\ref{sec:matrixForm},
before extensive numerical results are reported in Section~\ref{sec:numerics}.
We provide iteration counts for the Richardson and \textsc{GMRes} algorithms in
2D and 3D configurations with cross-points with physical boundaries as well as
interior cross-points. The influence of several impedance operators with respect
to different parameters: typical mesh size, wave number, number of subdomains
and varying coefficients (heterogeneous medium) is
studied.

\section{Problem under study}\label{PbUnderStudy}
We consider a very classical boundary value problem modeling scalar wave propagation
in an \textit{a priori}
heterogeneous medium in $\RR^{d}$ with $d=1,2$ or $3$. The computational domain $\Omega\subset \RR^{d}$
will be assumed bounded and polygonal ($d=2$) or polyhedral ($d=3$) for the sake of simplicity.
The material characteristics of the propagation medium will be represented by two functions
satisfying the following assumptions.

\begin{as}\label{Hypo1}\quad\\
  The functions $\kappa: \Omega\to \CC$ and $\mu:\Omega\to (0,+\infty)$ are measurable
  and satisfy 
  \begin{itemize}
  \item[(i)]   $\mrm{\Im m}\{\kappa(\bx)\}\geq 0, \mrm{\Re e}\{\kappa(\bx)\}\geq 0\;\forall \bx\in \Omega$
  \item[(ii)]  $\kappa_\infty:=\max(1,\Vert \kappa\Vert_{\mL^{\infty}(\Omega)}) <+\infty$
  \item[(iii)] $\Vert \mu \Vert_{\mL^{\infty}(\Omega)} + \Vert \mu^{-1}\Vert_{\mL^{\infty}(\Omega)}<+\infty$.
  \end{itemize}
\end{as}

\noindent 
These are both general and physically reasonable assumptions. Condition \textit{(i)} above
implies in particular that $\mrm{\Im m}\{\kappa^{2}(\bx)\}\geq 0$ and 
$\mrm{\Im m}\{\iu \kappa(\bx)\}\geq 0$ for all $\bx\in \Omega$. This means that 
the medium can only absorb or propagate energy. The above assumptions cover the case of
discontinuous material coefficients - piecewise constants for example. In addition, we consider
source terms $f\in \mL^{2}(\Omega)$ and $g\in \mL^{2}(\partial\Omega)$. The boundary value problem
under consideration will be 
\begin{equation}\label{InitialPb}
  \left\{\;\begin{aligned}
    & \text{Find\;$u\in \mH^{1}(\Omega)$ such that}\\
    & -\mrm{div}(\mu\nabla u) - \kappa^{2}u = f\quad \text{in}\;\Omega,\\
    & \left(\mu\partial_{\bn}-\iu\kappa\right)u = g\quad \text{on}\;\partial\Omega. 
  \end{aligned}\right.
\end{equation}
where \(\bn\) refers to the outward unit normal vector to \(\partial\Omega\).
As usual, for any domain $\omega\subset \Omega$, we consider the Sobolev space
$\mH^{1}(\omega):= \{v\in \mL^{2}(\omega): \nabla v\in\mL^{2}(\omega)\}$, and the space
of Dirichlet traces $\mH^{1/2}(\partial\omega):=\{ v\vert_{\partial\omega}:\;v\in\mH^{1}(\omega)\}$.
Because we are considering a wave propagation problem, inspired by \cite{zbMATH07248609}, we shall
equip those two spaces with the following norms 
\begin{equation}\label{H1Norm}
  \begin{aligned}
    & \Vert v\Vert_{\mH^{1}(\omega)}^{2}:=\Vert   \nabla v\Vert_{\mL^{2}(\omega)}^{2} +
    \kappa_\infty^2 \Vert v \Vert_{\mL^{2}(\omega)}^{2},\\
    & \Vert q\Vert_{\mH^{1/2}(\partial\omega)}:=\min \{\Vert v\Vert_{\mH^{1}(\omega)}:
    v\in \mH^1(\omega), v\vert_{\partial\omega} = q\}.
  \end{aligned}
\end{equation}
These norms involve some dependency with respect to the wave number $\kappa$ which
appears as an appropriate setting to derive frequency uniform results, see  e.g. \cite{zbMATH05969646}.
As usual, Problem~\eqref{InitialPb} can be put in variational form:
Find\;$u\in \mH^{1}(\Omega)$ such that $a(u,v) = \ell(v)$ $\forall v\in \mH^{1}(\Omega)$
where
\begin{equation}\label{SesquilinearForm}
  \begin{array}{rl}
      a(u,v):= & \int_{\Omega}\mu \nabla u\cdot
    \nabla \overline{v} -\kappa^{2} u\overline{v} \;\mathrm{d}\bx\;
        -\iu\int_{\partial\Omega} \kappa u \overline{v}\;\mathrm{d}\sigma \\[3pt]
        \ell(v):= &  \int_{\Omega}f\overline{v} \;\mathrm{d}\bx +
        \int_{\partial\Omega} g\overline{v}\;\mathrm{d}\sigma
  \end{array}  
\end{equation}
We are particularly interested in an effective numerical solution to~\eqref{InitialPb},
so we assume given a regular simplicial triangulation $\mathcal{T}_{h}(\Omega)$ of the domain
$\Omega$, and we assume $\overline{\Omega} = \cup_{\tau\in \mathcal{T}_{h}(\Omega)}\overline{\tau}$.
We shall denote $\Vh(\Omega)\subset \mH^1(\Omega)$ the space generated by conforming
$\mathbb{P}_{k}$-Lagrange functions (\(k\geq 1\)) constructed on
$\mathcal{T}_{h}(\Omega)$~\cite{Ern2004},
\begin{equation*}
  \Vh(\Omega):=\{v_h\in\mathscr{C}^0(\overline{\Omega}):\; v_h\vert_{\tau}\in
  \mathbb{P}_{k}(\tau),\forall \tau\in \mathcal{T}_{h}(\Omega)\}.
\end{equation*}
If $\omega\subset \Omega$ is any open subset that is resolved by the triangulation
i.e. $\overline{\omega} = \cup_{\tau\in \mathcal{T}_{h}(\omega)}\overline{\tau}$,
where $\mathcal{T}_{h}(\omega):=\{ \tau\in \mathcal{T}_{h}(\Omega): \tau\subset
\overline{\omega}\}$, then we denote
$\Vh(\omega):=\{ \varphi\vert_{\omega}:\;\varphi\in \Vh(\Omega)\}$
and $\Vh(\partial\omega):=\{ \varphi\vert_{\partial\omega}:\;\varphi\in \Vh(\Omega)\}$.
We will focus on the discrete variational formulation
\begin{equation}\label{DiscreteVF}
  \begin{array}{l}
    \text{Find}\;u_h \in\Vh(\Omega)\;\text{such that}\\
    a(u_h,v_h) = \ell(v_h)\quad \forall v_h\in\Vh(\Omega).
  \end{array}
\end{equation}
Devising efficient domain decomposition algorithms to solve this discrete problem
is the main goal of the present article. We make the general assumption that the material
characteristics and the mesh width are chosen so as to guarantee unique solvability of the
discrete problem.
\noindent
\begin{as}\label{Hypo2}\quad\\
$\alpha_{h}:= \displaystyle{\mathop{\textcolor{white}{p}\inf}_{u\in \Vh(\Omega)\setminus\{0\}}
 \sup_{v\in \Vh(\Omega)\setminus\{0\}}\frac{\vert a(u,v)\vert}
    {\Vert u\Vert_{\mH^{1}(\Omega)}\Vert v\Vert_{\mH^{1}(\Omega)}}}>0$.
\end{as}
This assumption simply means that the finite element matrix associated to \eqref{DiscreteVF}
is invertible. This is a reasonable assumption and hardly any numerical analysis seems possible without it. 
As a consequence, in the remaining of this article, we shall systematically assume $\alpha_h>0$.

\begin{rmrk}
In the general case, there is a mesh-threshold $\overline{h}>0$ such that $\alpha_h>0$ for $h\in (0,\overline{h})$.
The mesh-threshold is a function that depends on the material characteristics and the geometry
$\overline{h} = \overline{h}(\kappa,\mu,\Omega)$ so Assumption \ref{Hypo2} should be understood as a
condition imposed on $h,\mu,\kappa,\Omega$ all together; this is not just a condition on $h$ alone.
\end{rmrk}

\begin{rmrk}
  For fixed material characteristics $\mu,\kappa$ and fixed geometry $\Omega$, it is a classical
  consequence of Lax-Milgram's and Cea's lemma (see e.g.~\cite[chap.2]{MR1639879}) that unique
  solvability of Problem~\eqref{InitialPb} implies unique solvability of~\eqref{DiscreteVF}
  for a sufficiently fine mesh, and that the sesquilinear form $a(\cdot,\cdot)$ satisfies
  a uniform discrete lower bound $\alpha_{\star}:=\liminf_{h\to 0}\alpha_{h}>0$.
\end{rmrk}

\section{Geometric partitioning}\label{GeomPart}

We wish to describe and analyze a particular strategy for the solution
of Problem~\eqref{DiscreteVF} based on domain decomposition. As a consequence
we assume that the computational domain admits the decomposition
\begin{equation}\label{SubdomainPartition}
  \begin{array}{ll}
    \overline{\Omega} = \cup_{j=1}^{\mJ}\overline{\Omega}_{j}^{h},
    & \text{with}\quad\Omega_{j}^{h}\cap\Omega_{k}^{h} = \emptyset\quad
      \text{for}\;j\neq k\\[5pt]
    \Sigma := \cup_{j=1}^{\mJ}\Gamma_{j}^{h},
    & \text{where}\quad
        \Gamma_{j}^{h}:=\partial\Omega_{j}^{h},
  \end{array}
\end{equation}
where each $\Omega_{j}^{h}\subset \Omega$ is itself a polyhedral domain
that is exactly resolved by the triangulation.
Figure~\ref{Fig1} below gives
examples of the type of triangulation we consider.

A typical situation occurs when the computational domain is decomposed as
a first step in subdomains and, only afterwards, the mesh is generated in each subdomain
separately. In this case, mesh-conformity between subdomain triangulations has to be
enforced at interfaces. An example of this situation is represented in Figure~\ref{subfig1}.
Each subdomain $\Omega_{j}^{h}$, and the subdomain partition itself, then remains
unchanged as $h\to 0$. This is the case when the following condition is satisfied.

\begin{cond}[\textbf{fixed partition}]\label{Hypo4}\quad\\
    The subdomains $\Omega_{j}^{h} = \Omega_{j},\;j=1,\dots, \mJ$ are independent
    of the triangulation $\mathcal{T}_{h}(\Omega)$.
\end{cond}

\noindent 
However, we will \textit{not} consider that Condition~\ref{Hypo4} holds in
general (except in the examples of Section~\ref{sec:impedance} and in
Section~\ref{sec:fixedpartition}), because this partitioning approach is not
the most convenient from a practical viewpoint. This is the reason why we refer
to it as a ``condition'' instead of an ``assumption''.  

Another approach consists in generating a
mesh on the whole computational domain $\Omega$ first, and then subdividing it in subdomains by means of a graph
partitioner such as e.g. \textsc{Metis}~\cite{Karypis1999}.
In this manner, conformity of subdomain triangulations at interfaces is automatically
satisfied. However the partition itself has no reason to stabilize for $h\to 0$,
and there is no guarantee that the geometry of each subdomain converges. Boundaries
of subdomains may get rougher as $h\to 0$. This second situation is depicted in
Figure~\ref{subfig2}.
\begin{figure}[h]
  \hspace{0.15\textwidth}
  \begin{subfigure}{.3\textwidth}
    \hspace{0.1cm}
    \includegraphics[height=4cm]{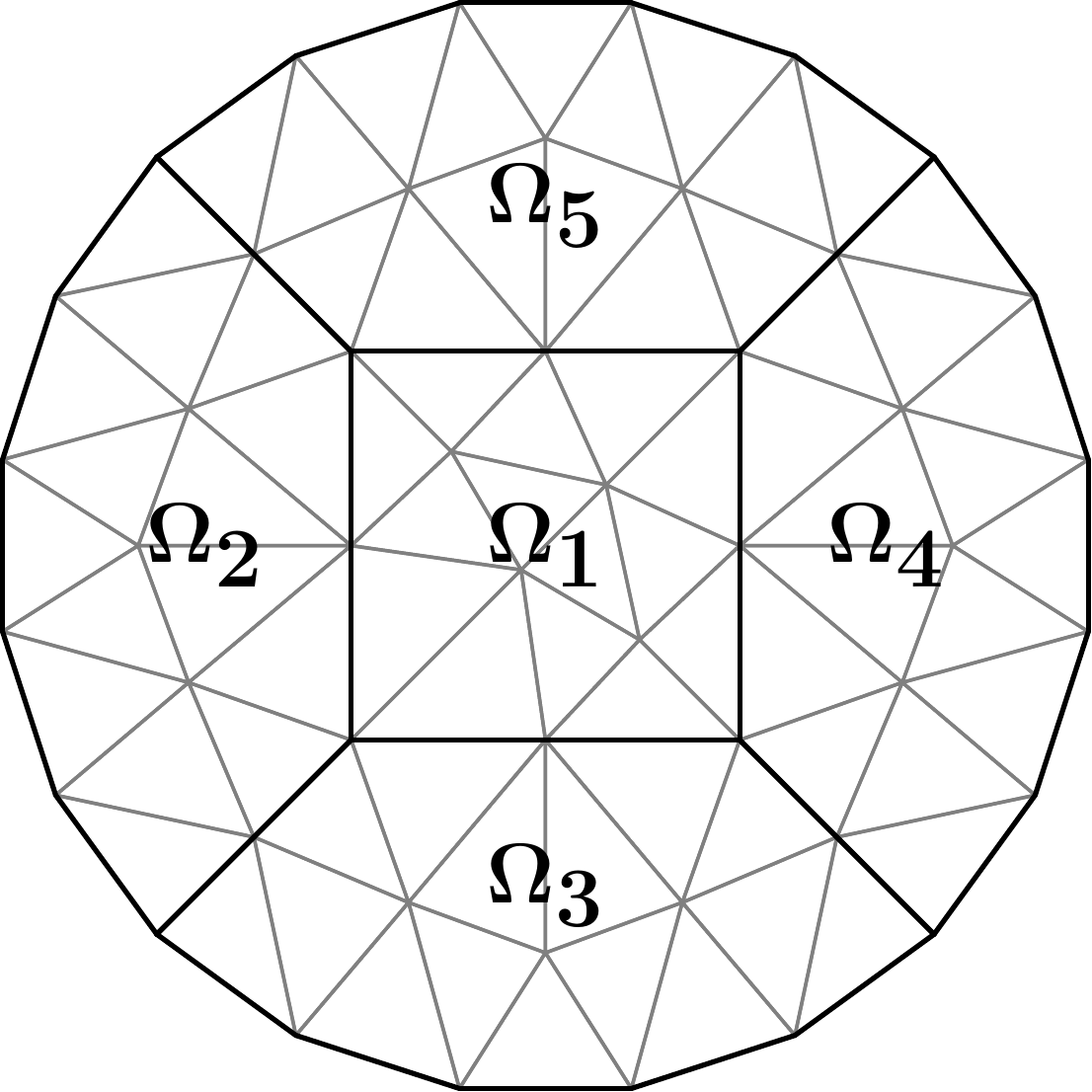}
    \caption{}\label{subfig1}
  \end{subfigure}
  \hspace{1cm}
  \begin{subfigure}{.3\textwidth}
    \hspace{0.1cm}
    \includegraphics[height=4cm]{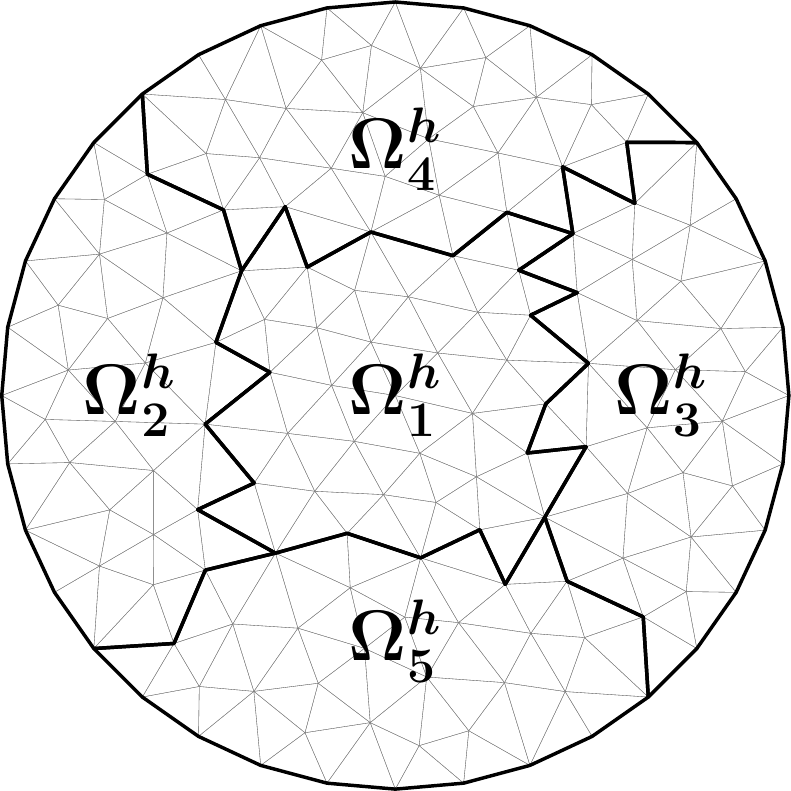}
    \caption{}\label{subfig2}
  \end{subfigure}
  \caption{Two approaches for partitioning the computational domain.}\label{Fig1}
\end{figure}

The analysis we present here covers both situations (a) and (b).
Our geometrical setting allows the presence of cross-points on the skeleton $\Sigma$ i.e.\ 
points where at least three subdomains may be adjacent. Because each $\Omega_{j}^{h}$ is
resolved by the triangulation, so is each boundary
$\Gamma_{j}^{h}$ as well as the skeleton $\Sigma$.

\section{Multi-domain functional setting}\label{MultiDomainFunctionalSetting}
We introduce continuous and discrete function spaces naturally associated
to this multi-domain setting. First we define $\mbH^{1}(\Omega) :=
\{u\in \mL^{2}(\Omega):\; u\vert_{\Omega_{j}^h} \in \mH^{1}(\Omega_{j}^h)\;
\forall j\}$ which contains the usual Sobolev space
$\mH^1(\Omega)\subset \mbH^{1}(\Omega)$. Next we consider the associated
finite element spaces
\begin{equation}\label{BrokenSpaces}
  \begin{aligned}
    & \mbVh(\Omega) := \{u\in\mL^{2}(\Omega):\;u\vert_{\Omega_{j}^{h}} \in \Vh(\Omega_{j}^{h})\;\forall j\}\;\subset\;\mbH^{1}(\Omega)\\
    & \text{so that}\;\; \Vh(\Omega)=\mH^{1}(\Omega)\cap \mbVh(\Omega)\subset\mbVh(\Omega).
  \end{aligned}
\end{equation}
The space $\mbVh(\Omega)$ consists in functions that are $\mathbb{P}_1$-Lagrange piecewise
with respect to the subdomain partition. It is naturally isomorphic to $\Vh(\Omega_1)\times
\dots\times \Vh(\Omega_\mJ)$. The elements of $\mbVh(\Omega)$ may possibly admit
Dirichlet jumps across interfaces between subdomains, while $\Vh(\Omega)$ is a closed subspace of
$\mbVh(\Omega)$ characterized by the constraint that Dirichlet traces match across
interfaces $\Gamma_{j}^{h}\cap\Gamma_{k}^{h}$.  Since we will be interested in transmission
conditions, we need to consider the trace operation
\begin{equation*}
  v\vert_{\Sigma}:= (v\vert_{\Gamma_{1}^{h}},\dots,v\vert_{\Gamma_{\mJ}^{h}})
\end{equation*}
where each $v\vert_{\Gamma_{j}^h}$ is a Dirichlet trace taken \textit{from the interior}
of the subdomain $\Omega_{j}^{h}$, which is a convention we will systematically adopt
through this article. The trace operator above continuously maps $\mbVh(\Omega)$
(resp. $\Vh(\Omega)$) onto the space of Dirichlet multi-traces (resp. single-traces)
\begin{equation}\label{MultiTraceSpaces}
  \begin{aligned}
    \mbVh(\Sigma)
    & := \{v\vert_{\Sigma}:\; v\in  \mbVh(\Omega)\} = 
    \Vh(\Gamma_{1}^{h})\times\cdots\times \Vh(\Gamma_{\mJ}^{h})\\
    \Vh(\Sigma)
    & := \{v\vert_{\Sigma}:\; v\in  \Vh(\Omega)\}\subset\mbVh(\Sigma).
  \end{aligned}
\end{equation}
where $\Vh(\Gamma_{j}^{h})=\{ v\vert_{\Gamma_{j}^{h}}: v\in \Vh(\Omega)\}$, see Section \ref{PbUnderStudy}.
In notations \eqref{BrokenSpaces}-\eqref{MultiTraceSpaces}, we used the letter "$\mbV$"
for reference to spaces where no continuity constraint is imposed at interfaces, and we used
the letter "$\mV$" for reference to spaces where continuity holds through interfaces.
The spaces $\mbVh(\Omega)$ and $\mbVh(\Sigma)$ are naturally equipped with cartesian product norms
stemming from \eqref{H1Norm}: for $ u\in \mbVh(\Omega)$ and
$\ctrq = (q_1,\dots,q_\mJ)\in \mbVh(\Sigma)$ we set
\begin{equation}\label{ContinuousMultiDomainNorms}
  \begin{aligned}
    & \Vert u\Vert_{\mbH^{1}(\Omega)}^{2} := \Vert u|_{\Omega_{1}^{h}}\Vert_{\mH^{1}(\Omega_{1}^{h})}^{2}+
    \dots+ \Vert u|_{\Omega_{J}^{h}} \Vert_{\mH^{1}(\Omega_{\mJ}^{h})}^{2},\\
    & \Vert \ctrq\Vert_{\mbH^{1/2}(\Sigma)}^{2} := \Vert q_1\Vert_{\mH^{1/2}(\Gamma_{1}^{h})}^{2}+
    \dots+ \Vert q_\mJ \Vert_{\mH^{1/2}(\Gamma_{\mJ}^{h})}^{2}.
  \end{aligned}
\end{equation}

\section{Discrete harmonic extensions}\label{sec:tracenorms}

Besides \eqref{ContinuousMultiDomainNorms}, several choices of norms are possible for the
multi-trace space $\mbVh(\Sigma) = \Vh(\Gamma_1^h)\times \dots\times \Vh(\Gamma_\mJ^h)$.
Because we are dealing with a discrete variational formulation, there is another norm,
of a more discrete nature, that will arise as pivotal in our analysis. To describe
it, we first introduce norm minimizing lifting operators $\rho_j:\Vh(\Gamma_j^h)\to \Vh(\Omega_j^h)$
defined, for all $q\in \Vh(\Gamma_j^h)$, by
\begin{equation}\label{HarmonicExtensionMap}
  \begin{aligned}
    & \rho_j(q)\in \Vh(\Omega_j^h),\;\rho_j(q)\vert_{\Gamma_j^h} = q\;\;\text{and}\\
    & \Vert \rho_j(q)\Vert_{\mH^1(\Omega_j^h)} = \min\{ \Vert v\Vert_{\mH^1(\Omega_j^h)}:\; v\in \Vh(\Omega_j^h), v\vert_{\Gamma_j^h} = q\}. 
  \end{aligned}
\end{equation}
Such operators are commonly encountered in domain decomposition literature
(see e.g. \cite[\S 4.4]{MR2104179} or \cite[1.2.6]{MR3013465}) and are sometimes
referred to as harmonic extension maps. These can be collected into a global
lifting operator defined by
\begin{equation*}
  \begin{aligned}
    & \rho_h:\mbVh(\Sigma)\to \mbVh(\Omega),\\
    & \rho_h(\ctrq)\vert_{\Omega_j} := \rho_j(q_j)\;\;\forall j=1\dots \mJ 
  \end{aligned}
\end{equation*}
for all $\ctrq = (q_1,\dots,q_\mJ)$. By construction we have
$\rho_h(\ctrq)\vert_{\Sigma} = \ctrq$ and, in particular,
$\rho_h$ is one-to-one. A possible norm for discrete Dirichlet
traces stems directly from discrete harmonic
extension maps. For $\ctrq = (q_1,\dots,q_\mJ)\in \mbVh(\Sigma)$ we set
\begin{equation}\label{MTFNorm}
  \Vert \ctrq\Vert_{\rh}^2 := \Vert \rho_h(\ctrq)\Vert_{\mbH^1(\Omega)}^2 :=
  \sum_{j=1}^{\mJ}\Vert \rho_j(q_j)\Vert_{\mH^1(\Omega_j^h)}^2.
\end{equation}
By the very definition of this norm, the operator $\rho_h$ has unitary continuity modulus.
The norm \eqref{MTFNorm} inherits wave number dependency from \eqref{H1Norm}. It also a priori
depends on the triangulation because minimization in \eqref{HarmonicExtensionMap} is performed
over the discrete space $\Vh(\Omega_j^h)$.

\quad\\
Comparing \eqref{HarmonicExtensionMap}-\eqref{MTFNorm} on one hand, and
\eqref{H1Norm}-\eqref{ContinuousMultiDomainNorms} on the other hand,
the norm $\Vert \cdot\Vert_{\rh}$ appears as a discrete counterpart of
$\Vert \cdot\Vert_{\mbH^{1/2}(\Sigma)}$. For fixed material characteristics these
two norms are actually $h$-uniformly equivalent.

\begin{lem}\label{ScottZhangRobustness}\quad\\
  The continuity of the trace operator systematically implies $\Vert \ctrq\Vert_{\mbH^{1/2}(\Sigma)}\leq \Vert \ctrq\Vert_{\rh}$ for all $\ctrq\in \mbVh(\Sigma)$.
  Assume in addition that the family of triangulations $\{\mathcal{T}_{h}(\Omega)\}_{h>0}$ is $h$-uniformly shape regular,
  \begin{equation}\label{ShapeRegularity}
    \liminf_{h\to 0,\tau \in \mathcal{T}_{h}(\Omega)}
    \sup\{ \mrm{diam}(\mrm{B}):\;\mrm{B}\;\text{is a ball}\,\subset \tau\}/\mrm{diam}(\tau)>0
  \end{equation}
  and that the material characteristics (in particular the wave number) are fixed as well as the computational
  domain $\Omega$. Then the discrete harmonic extension operator is  $h$-uniformly bounded
  \begin{equation*}
    \limsup_{h\to 0}\sup_{\ctrq\in \mbVh(\Sigma)\setminus\{0\}}\frac{\hspace{-0.55cm}\Vert \ctrq\Vert_{\rh}}{\Vert \ctrq\Vert_{\mbH^{\frac{1}{2}}(\Sigma)}} <+\infty.
  \end{equation*}
\end{lem}
\noindent \textbf{Proof:}

For any $p\in \mH^{1/2}(\Gamma_j^h)$ define $\theta_j(p)$ as the unique element of 
$\mH^{1}(\Omega_j^h)$ satisfying $\theta_j(p)\vert_{\Gamma_j^h} = p$ and
$\Vert \theta_j(p)\Vert_{\mH^1(\Omega_j^h)} = \min\{\Vert v\Vert_{\mH^1(\Omega_j^h)}:\, v\in \mH^1(\Omega_j^h),
\,v\vert_{\Gamma_j^h} = p\}$, and then define $\theta:\mbH^{1/2}(\Sigma)\to \mbH^1(\Omega)$ by
$\theta(q_1,\dots,q_\mJ)\vert_{\Omega_j^h} := \theta_j(q_j)$ for all $j=1\dots \mJ$.
Here the operator \(\theta\) is similar to the operator \(\rho_{h}\)
in the sense that it is also an harmonic extension map, but defined on the continuous spaces.
By definition of the
norms \eqref{H1Norm}-\eqref{ContinuousMultiDomainNorms} we have
for any \(\ctrq\in\mbH^{1/2}(\Sigma)\)
\begin{equation}
\Vert \ctrq\Vert_{\mbH^{1/2}(\Sigma)} = \Vert \theta(\ctrq)\Vert_{\mbH^1(\Omega)}.
\end{equation}
Next let $\pi_{h,j}:\mH^{1}(\Omega_{j}^{h})\to \Vh(\Omega_{j}^{h})$ refer
to the Scott-Zhang interpolation  operator~\cite{zbMATH04141375} and~\cite[\S4.8]{zbMATH05223061}. This operator
satisfies $\pi_{h,j}(v) = v$ for all $v\in \Vh(\Omega_{j}^{h})$,  and it is uniformly bounded
\begin{equation}\label{hUniformContinuity}
  \limsup_{h\to 0}\sup_{v\in\mH^{1}(\Omega_{j}^{h})\setminus\{0\}}
  \Vert \pi_{h,j}(v)\Vert_{\mH^{1}(\Omega_{j}^{h})}/ \Vert v\Vert_{\mH^{1}(\Omega_{j}^{h})}<+\infty.
\end{equation}
This modulus of the Scott-Zhang interpolator only depends on the shape regularity constant~\eqref{ShapeRegularity}
of the triangulation and is thus bounded independently
of the shape of the subdomains. Besides it can be defined in such a way that it guarantees
$\pi_{j,h}(v)\vert_{\Gamma_{j}^{h}} = v\vert_{\Gamma_{j}^{h}}$ if $v\vert_{\Gamma_{j}^{h}}\in \Vh(\Gamma_{j}^{h})$.
Define $\pi_h:\mbH^1(\Omega)\to\mbVh(\Omega)$ by collecting these local Scott-Zhang operators
$\pi_h(v)\vert_{\Omega_j^h}:= \pi_{j,h}(v\vert_{\Omega_j^h})$ for all $j=1\dots\mJ$. From 
\eqref{hUniformContinuity} we deduce
\begin{equation}
  C_\pi= \limsup_{h\to 0}\sup_{v\in\mbH^{1}(\Omega)\setminus\{0\}}
  \Vert \pi_{h}(v)\Vert_{\mbH^{1}(\Omega)}/ \Vert v\Vert_{\mbH^{1}(\Omega)}<+\infty.
\end{equation}
Finally, taking  account of the definition of $\Vert \cdot\Vert_{\rh}$ given by
\eqref{HarmonicExtensionMap}-\eqref{MTFNorm},  we obtain
$\Vert \ctrq\Vert_{\rh}= \Vert \rho_h(\ctrq)\Vert_{\mbH^1(\Omega)}
\leq \Vert \pi_{h}(\theta(\ctrq))\Vert_{\mbH^{1}(\Omega)} \leq C_\pi \Vert \theta(\ctrq)\Vert_{\mbH^{1}(\Omega)}
= C_\pi\Vert \ctrq\Vert_{\mbH^{1/2}(\Sigma)}$
for any \(\ctrq\in\mbH^{1/2}(\Sigma)\).\qed

\quad\\
To conclude this paragraph, let us point that the constant associated to the $\limsup_{h\to0}$ in Lemma \ref{ScottZhangRobustness}
above a priori depends on the material characteristics and in particular the wave number.

\section{Impedance operator}\label{sec:impedance}
To define a scalar product on the space of traces, one may of course consider the scalar
product associated with \eqref{MTFNorm} or \eqref{ContinuousMultiDomainNorms}. Many other
choices are possible that appear more convenient in a finite element context. This is why we
will \textit{a priori} consider another one.
\begin{as}\label{Hypo5}\quad\\
  The sesquilinear form $t_{h}(\cdot,\cdot):\mbVh(\Sigma)\times \mbVh(\Sigma)\to \CC$ is a scalar product,
  and the associated norm is denoted
  \begin{equation*}
    \Vert \ctrw\Vert_{\Th}:= \sqrt{t_{h}(\ctrw,\ctrw)}.
\end{equation*}
\end{as}
We state this as an assumption because the sesquilinear form $t_h(\cdot,\cdot)$, that we shall later
refer to as the \textit{impedance}, will play a central role in our theory, and it is important to underline 
that our analysis allows rather general form for the impedance.

We insist that no other assumption  on the impedance is needed for the subsequent analysis.
In particular we do not assume that $t_{h}(\cdot,\cdot)$ satisfies a uniform discrete inf-sup
condition so that $\Vert \cdot\Vert_{\Th}$ need not be $h$-uniformly equivalent to
$\Vert \cdot\Vert_{\mbH^{1/2}(\Sigma)}$ or $\Vert \cdot\Vert_{\rh}$. The following quantities
will come into play
\begin{equation}\label{StabilityContinuityModulus}
  \lambda_{h}^{-}:=\inf_{\ctrw\in\mbVh(\Sigma)\setminus\{0\}}\frac{\Vert \ctrw\Vert_{\Th}}{\Vert \ctrw\Vert_{\rh}}
  \qquad\text{and}\qquad
  \lambda_{h}^{+}:=\sup_{\ctrw\in\mbVh(\Sigma)\setminus\{0\}}\frac{\Vert \ctrw\Vert_{\Th}}{\Vert \ctrw\Vert_{\rh}}.
\end{equation}
Because $t_{h}(\cdot,\cdot)$ is supposed to be a scalar product, we have $0<\lambda_{h}^{-}\leq \lambda_{h}^{+}<+\infty$
for each $h$. However, in certain cases of practical importance, one may have $\liminf_{h\to 0}\lambda_{h}^{-} = 0$ or
 $\limsup_{h\to 0}\lambda_{h}^{+} = +\infty$. To fix ideas, we give a few examples of
possible impedance operators.

\begin{exmpl}[\textbf{Despr{\'e}s impedance}]\label{despresimpedance}
Consider a fixed parameter $\kappa_{\rref}>0$ that will serve as a reference
wave number. A first possible choice consists in defining impedance using
surface mass matrices as follows
\begin{equation}\label{DespresOperator}
    t_{h}(\ctrp,\ctrq) = \sum_{j=1}^{\mJ}\int_{\Gamma_{j}}\kappa_{\rref}\,p_{j}\overline{q}_{j} \;\mathrm{d}\sigma.
\end{equation}
This choice of impedance was the one originally introduced in~\cite{MR1071633,MR1105979,MR1291197,MR1227838}. 

It is a scalar product on $\mbVh(\Sigma)$ that is $h$-uniformly bounded $\limsup_{h\to 0}\lambda_{h}^{+}<+\infty$. 
If Condition~\ref{Hypo4} holds and material characteristics are fixed, it does not satisfy any $h$-uniform
discrete inf-sup condition. Inverse estimates show that there exists a constant $C>0$ independent of $h$ such
that  $\lambda_{h}^{-}\geq C \sqrt{h}$, see e.g~\cite[Lemma 4.11]{MR2104179}. 
\end{exmpl}

%\medskip
\begin{exmpl}[\textbf{Second order differential operator}]\label{sndorderimpedance}
Consider two constants $a,b>0$ that, in practice, are fitting parameters requiring
calibration.  Another choice of impedance is based on an order 2 surface differential operator
and involves both mass and stiffness matrices
\begin{equation}\label{ScndOrder}
  t_{h}(\ctrp,\ctrq) = \sum_{j=1}^{\mJ}\int_{\Gamma_{j}}
  a \nabla_{\Gamma_{j}} p_{j}\cdot \nabla_{\Gamma_{j}}\overline{q}_{j} + b\,
  p_{j}\overline{q}_{j} \;\mathrm{d}\sigma.
\end{equation}
In the choice above, the operators $\nabla_{\Gamma_{j}}$ refer to surface gradients.
Many variants of this condition can be considered. The coefficients $a,b$ may vary from one subdomain
to another i.e. $a = a_{j}$ (resp. $b = b_{j}$) on $\Gamma_{j}$. Also these coefficients may depend on
the mesh width $h$ and the wave number $\kappa$. Such a choice of impedance (or a variant of it) was
considered e.g.\ in~\cite{MR1924414,MR2097757,MR3549886}.

Impedance~\eqref{ScndOrder} also yields a scalar product on $\mbVh(\Sigma)$. If Condition~\ref{Hypo4} holds
and material characteristics are fixed, it can be proved to be $h$-uniformly inf-sup stable $\liminf_{h\to 0}\lambda_{h}^{-}>0$,
but it is not $h$-uniformly bounded anymore. Once again inverse estimates yield the existence of $C>0$ independent of
$h$ such that $\lambda_{h}^{+}\leq C/\sqrt{h}$.
\end{exmpl}

%\medskip
\begin{exmpl}[\textbf{Integral operator based impedance}]\label{hypersingularimpedance}
Another possibility is an impedance based on some integral operator. Consider parameters
$a,\delta>0$. The parameter $\delta$ would represent how localized the kernel is. To mimic the properties
of the trace norm, one may consider an analytical expression based on the
hypersingular operator see e.g.~\cite[\S 3.3.4]{MR2743235} or~\cite[\S
6.5]{MR2361676}. For problems posed in $\RR^{3}$ this would correspond to 
\begin{equation}\label{HyperSingular}
  \begin{aligned}
    & t_{h}(\ctrp,\ctrq) = \sum_{j=1}^{\mJ} \int_{\Gamma_{j}\times\Gamma_{j}}
    a\;\frac{\exp(-\vert \bx-\by\vert/\delta)}{4\pi \vert \bx-\by\vert}
    \big(\; \mrm{curl}_{\Gamma_{j}}p_j(\bx)\cdot
    \mrm{curl}_{\Gamma_{j}}\overline{q}_j(\by)\\
    &  \hspace{5cm}+ \delta^{-2} \bn_{j}(\bx)\cdot\bn_{j}(\by) \,p_j(\bx)\overline{q}_j(\by) \;\big)\;\mathrm{d}\sigma(\bx,\by),
  \end{aligned}
\end{equation}
where $\bn_{j}(\bx)$ is the vector normal  to $\Gamma_{j}$ directed toward
the exterior of $\Omega_{j}$, and $\mrm{curl}_{\Gamma_{j}}v(\bx) := \bn_{j}(\bx)\times
\nabla_{\Gamma_{j}}v(\bx)$ is the surface curl. Such a choice was considered
in~\cite{LecouvezCollinoJoly,lecouvez:tel-01444540}, and  a variant of it was proposed
in~\cite{MR3989867} in the context of Maxwell's equations. As discussed in
\cite[\S 3.1.3]{lecouvez:tel-01444540}, the motivation for considering
non-local impedances such as \eqref{HyperSingular} is that this leads to $h$-uniform
rate of convergence of the OSM algorithm, in contrast with local impedances such as
\eqref{DespresOperator} or \eqref{ScndOrder}.  

Once again~\eqref{HyperSingular} is a scalar product on $\mbVh(\Sigma)$. Because~\eqref{HyperSingular} stems
from an operator that is both bounded and coercive in the continuous trace space, if Condition~\ref{Hypo4} holds
and material characteristics are fixed, then  $0<\liminf_{h\to 0}\lambda_{h}^{-}\leq \limsup_{h\to 0}\lambda_{h}^{+}<+\infty$. 
\end{exmpl}

\begin{exmpl}[\textbf{Schur complement based impedance}]\label{ex:schurcomplementimpedance}
Another non-local operator is given by a Schur complement based impedance
operator, which corresponds to the scalar product associated to the following norm:
\begin{equation}\label{SchurComplementBasedImpedance}
    t_h(\ctrp,\ctrp) =
    \Vert \ctrp\Vert_{\rh}^2 =
    \Vert \rho_h(\ctrp)\Vert_{\mbH^1(\Omega)}^2 =
    \sum_{j=1}^\mJ\Vert \rho_j(p_j)\Vert_{\mH^1(\Omega_j^h)}^2,
\end{equation}
where $\ctrp = (p_{1},\dots, p_{\mJ})\in\mbVh(\Sigma)$.
In practice, computing this scalar product involves computing the local
solutions \(\rho_{j}(p_j)\), see \eqref{HarmonicExtensionMap},  which can be
done in parallel using the same numerical scheme and mesh as for the
solutions \(u_{j}\) of the propagative local sub-problems.

From the definition of $\lambda_h^\pm$ given by~\eqref{StabilityContinuityModulus}
we readily obtain $\lambda_h^\pm = 1$ in this case, which implies unconditional
\(h\)-uniform stability for this choice of operator.
\end{exmpl}

\noindent 
In all the examples we gave above, the impedance does not couple distinct subdomains i.e.\ 
it takes the form of a sum of local contributions, which is rather
natural in the context of domain decomposition.

\begin{defn}[\textbf{Diagonal impedance}]\label{DiagoImpedance}\quad\\
  A scalar product $t_{h}(\cdot,\cdot):\mbVh(\Sigma)\times \mbVh(\Sigma)\to \CC$ will
  be called diagonal if there are local scalar products  $t_{h}^{j}(\cdot,\cdot):\Vh(\Gamma_{j}^{h})
  \times\Vh(\Gamma_{j}^{h})\to \CC$ such that, for $\ctrp = (p_{1},\dots, p_{\mJ})$ and
  $\ctrq = (q_{1},\dots, q_{\mJ})$, we have  
  \begin{equation*}
    t_{h}(\ctrp,\ctrq) = t_{h}^{1}(p_{1},q_{1})+\cdots + t_{h}^{\mJ}(p_{\mJ},q_{\mJ}).
  \end{equation*}
\end{defn}

\noindent 
Many other choices of impedance operator are possible, see e.g.~\cite{MR1764190,lecouvez:tel-01444540}. 
We do not assume a priori that the impedance $t_h(\cdot,\cdot)$ is diagonal because this
is not needed for the subsequent analysis. Later on though, we shall point the practical benefit of
considering diagonal impedances, see Remark \ref{BenefitDiagonality}.

\quad\\
To conclude this section we point out that the scalar product $t_{h}(\cdot,\cdot)$ will be used to represent
linear functionals associated to boundary terms.  

\begin{lem}\label{LiftingLemma}\quad\\
  Under Assumption \ref{Hypo5}, assume that $\phi:\mbVh(\Omega)\to \CC$
  is a linear form satisfying $\phi(w) = 0$ for all $w\in \mbVh(\Omega)$ such that
  $w\vert_{\Sigma} = 0$. Then there exists a unique $\ctrq\in \mbVh(\Sigma)$ such that
  $\phi(v) = t_{h}(v\vert_{\Sigma},\ctrq)\forall v\in \mbVh(\Omega)$.
\end{lem}
\noindent \textbf{Proof:}

Since $w=\rho_{h}(v\vert_{\Sigma}) - v\in \mbVh(\Omega)$ satisfies $w\vert_{\Sigma} = 0$ for all $v\in\mbVh(\Omega)$,
we deduce that $\phi(v) = \phi(\rho_{h}(v\vert_{\Sigma}))$ for all $v\in \mbVh(\Omega)$. Since $t_{h}(\cdot,\cdot)$
is a scalar product on $\mbVh(\Sigma)$, there exists a unique $\ctrq\in \mbVh(\Sigma)$ such that
$t_{h}(\ctrp,\ctrq) = \phi(\rho_{h}(\ctrp))\forall \ctrp\in\mbVh(\Sigma)$ by Riesz representation. 
\qed%

\section{Characterization of interface conditions}\label{sec:reformulation}

Transmission conditions are critical in domain decomposition because this is what induces
coupling between subdomains. In the present section, we will work out a new way to impose
them.

In \eqref{DiscreteVF}, Dirichlet transmission conditions are enforced by assuming
that traces match across interfaces. Such conditions are encoded into the single-trace space
introduced in \eqref{MultiTraceSpaces}, 
\begin{equation*}
    \Vh(\Sigma)=\{ {(u\vert_{\Gamma_{j}^{h}})}_{j=1}^{\mJ}
    \in \mbVh(\Sigma):\; u\in \Vh(\Omega)\}.
\end{equation*}
This space is a discrete counterpart of what has been referred to as "Dirichlet 
single-trace" space in the literature dedicated to Multi-Trace formalism, see~\cite{C11_144,MR3403719,MR3069956}. 
We first state a simple result that helps characterize $\Vh(\Omega)$ as a subspace of
$\mbVh(\Omega)$ by means of its traces on the skeleton.
\begin{lem}\label{TraceMatching}\quad\\
  Any $u\in\mbVh(\Omega)$ belongs to $\Vh(\Omega)$ if and only if
  $u\vert_{\Sigma}\in \Vh(\Sigma)$. 
\end{lem}

\noindent 
An easy consequence of the previous lemma is that $\rho_{h}(\ctrq)\in \Vh(\Omega)$ whenever
$\ctrq\in \Vh(\Sigma)$. Because we are aiming at domain decomposition, we wish to decouple
subdomains as much as possible, and thus eliminate any reference to $\Vh(\Sigma)$.
This is our motivation in searching for characterizations of
this space. The next result provides such a characterization by means of an ``exchange operator''
$\Pi$. Its proof is routine verification left to the reader. 

\begin{lem}[\textbf{Exchange operator}]\label{DefOrthoProj}\quad\\
  Under Assumption \ref{Hypo5}, define $\Pi:\mbVh(\Sigma)\to\mbVh(\Sigma)$
  such that $(\Id+\Pi)/2$ is the orthogonal projector onto $\Vh(\Sigma)$ for the scalar product
  \(t_{h}\). Then \(\Pi\) is an isometry, we have $\Pi^{2} = \Id$ and $\Vert
  \Pi(\ctru)\Vert_{\Th} = \Vert \ctru\Vert_{\Th}$ for all
  $\ctru\in \mbVh(\Sigma)$. Moreover for any pair $(\ctru,\ctrq)\in \mbVh(\Sigma)
  \times\mbVh(\Sigma)$ we have
  \begin{equation}\label{CaracTranCond}
    (\ctru,\ctrq)\in \Vh(\Sigma)
    \times\Vh{(\Sigma)}^{\perp}\quad\iff\quad -\ctrq+\iu\ctru = \Pi(\ctrq+\iu\ctru).
  \end{equation}
\end{lem}

\noindent 
Computing the action of the exchange operator through the operation $\ctrp\mapsto \Pi(\ctrp)$
is non-trivial from an effective computational viewpoint. This can be achieved through
an orthogonal projection onto the subspace $\Vh(\Sigma)$ which, variationally, rewrites as
follows
\begin{equation}\label{OrthoProj}
  \begin{aligned}
    & \quad\Pi(\ctrq) := -\ctrq + 2 \ctrp\quad \text{where}\;\;\ctrp\in \Vh(\Sigma)\;\;\text{solves}\\
    & \quad t_{h}( \ctrp ,\ctrw) =  t_{h}(\ctrq,\ctrw)\quad\forall\ctrw\in \Vh(\Sigma).
  \end{aligned}
\end{equation}
Of course this orthogonal projection requires solving a problem posed globally on the whole 
skeleton $\Sigma$. Here the choice of the scalar product $t_{h}$ does matter: it should
be chosen so that the orthogonal projection in~\eqref{OrthoProj} is easy to compute.
This variational problem makes the operator $\Pi$ a priori non-local. For certain choices 
of impedance, this exchange operator may couple distant subdomains that are not a priori adjacent. 
This will be a salient feature of our strategy, and a key difference in comparison with existing 
literature.

\begin{rmrk}\label{remark1}
  Admittedly, this non-locality raises a computational difficulty. However the
  variational problem~\eqref{OrthoProj} is symmetric positive definite, and
  takes the very same form as the Schur complement systems encountered in the
  analysis of substructuring methods, see~\cite[\S 4.3]{MR2104179},~\cite[\S
  2.1]{MR3013465} or~\cite[\S 6.4]{MR3450068}. Current literature offers very
  efficient scalable two level DDM preconditioners for tackling such a problem as 
  \eqref{OrthoProj}, like e.g. the GenEO approach \cite{zbMATH06669320}.
\end{rmrk}

\section{Reformulation of the wave propagation problem}

In the previous section we discussed in depth the effect of domain partitioning
on our discrete functional setting. This led us to propose, in Section
\ref{sec:reformulation}, a new way of imposing transmission conditions, see \eqref{CaracTranCond}.
We now apply the outcomes of this discussion to the discrete problem \eqref{DiscreteVF}
under study. The sesquilinear form $a(\cdot,\cdot)$ introduced in \eqref{SesquilinearForm}
extends as a map $a(\cdot,\cdot):\mbVh(\Omega)\times \mbVh(\Omega)\to \CC$ defined by
\begin{equation*}
  \begin{array}{rl}
    a(u,v) & := \sum_{j=1}^{\mJ}a_{\Omega_{j}^{h}}(u,v)\quad \text{where}\\[5pt]
    a_{\Omega_{j}^{h}}(u,v) & := \int_{\Omega_{j}^{h}}\mu\nabla u\cdot\nabla \overline{v}
    -\kappa^{2} u\overline{v} \;\mathrm{d}\bx
    -\iu\int_{\partial\Omega_{j}^{h}\cap\partial\Omega}\kappa u \overline{v}\;\mathrm{d}\sigma.
  \end{array}
\end{equation*}
Note that Assumptions \ref{Hypo1} imply $\Im m\{a(u,u)\}\leq 0\;\forall u\in \mbVh(\Omega)$.
Besides the inf-sup constant $\alpha_h$ introduced in Assumption \ref{Hypo2}, the material
characteristics $\mu,\kappa$ and the geometry of $\Omega$ enter our analysis
through the continuity modulus
\begin{equation}\label{ContinuityModulusA}
  \Vert a\Vert:=\sup_{u,v\in \mbH^{1}(\Omega)\setminus\{0\}}
  \frac{\vert a(u,v)\vert}{\Vert u\Vert_{\mbH^{1}(\Omega)}\Vert v\Vert_{\mbH^{1}(\Omega)}}.
\end{equation}

\begin{rmrk}\label{UniformBoundednessSesquilinearform}
The same arguments as in \cite[Lem.2.4]{zbMATH07248609}, and
in particular the multiplicative trace inequality from
\cite[last equation on p.41]{zbMATH05960425}, show that $\Vert a\Vert$
is bounded uniformly with respect to $\sup_\Omega \vert \kappa\vert$. 
\end{rmrk}
Like for $a(\cdot,\cdot)$, the functional $\ell(\cdot)$ induces a continuous map on $\mbVh(\Omega)$
defined for all  $v\in\mbVh(\Omega)$ by
\begin{equation*}
  \begin{array}{rl}
    \ell(v) & :=  \sum_{j=1}^{\mJ}\ell_{\Omega_{j}^{h}}(v)\quad \text{where}\\[5pt]
    \ell_{\Omega_{j}^{h}}(v) & := \int_{\Omega_{j}^{h}}f \overline{v} \;\mathrm{d}\bx +
    \int_{\partial \Omega_{j}^{h}\cap\partial \Omega }g\overline{v} \;\mathrm{d}\sigma.
  \end{array}
\end{equation*}
The previous notations combined with the impedance
$t_h(\cdot,\cdot)$ introduced in Section \ref{sec:impedance} can be used to
reformulate~\eqref{DiscreteVF} as a saddle-point problem involving the space
$\mbVh(\Omega)$ instead of $\Vh(\Omega)$, and enforcing transmission conditions by
means of Lagrange multipliers.

\begin{prop}\label{EquivalenceFormulation}\quad\\
  Let Assumption \ref{Hypo5} hold, and set ${\Vh(\Sigma)}^{\perp}:=\{ \ctru\in\mbVh(\Sigma):
  \;t_{h}(\ctru,\ctrv) = 0\; \forall \ctrv\in\Vh(\Sigma)\}$. If the function $u\in \Vh(\Omega)$
  solves~\eqref{DiscreteVF}, then there exists $\ctrq\in \Vh{(\Sigma)}^{\perp}$ such that
  \begin{equation}\label{SystemIntermediaire}
    \begin{aligned}
      & (u,\ctrq)\in \mbVh(\Omega)\times\Vh{(\Sigma)}^{\perp}\quad\text{and}\\
      & a(u,v) - t_{h}(\ctrq,v\vert_{\Sigma}) = \ell(v)\quad \forall v\in \mbVh(\Omega),\\
      & t_{h}(u\vert_{\Sigma},\ctrw) = 0\hspace{2.3cm} \forall \ctrw\in \Vh{(\Sigma)}^{\perp}.
    \end{aligned}  
  \end{equation}
  Reciprocally if $(u,\ctrq)\in \mbVh(\Omega)\times {\Vh(\Sigma)}^{\perp}$ solves~\eqref{SystemIntermediaire},
  then $u\in \Vh(\Omega)$ and it is the unique solution to Problem~\eqref{DiscreteVF}.
\end{prop}
\noindent \textbf{Proof:}

Assume first that $u\in\Vh(\Omega)$ solves~\eqref{DiscreteVF}. By definition
we have $u\vert_{\Sigma}\in \Vh(\Sigma)$ so that $t_{h}(u\vert_{\Sigma},\ctrw) = 0
\;\forall \ctrw\in {\Vh(\Sigma)}^{\perp}$. In addition,  applying Lemma~\ref{LiftingLemma}
to the functional $\phi(v):=a(u,v)-\ell(v), v\in \mbVh(\Omega)$,
we deduce that there exists $\ctrq\in \mbVh(\Sigma)$ such that $a(u,v)-\ell(v)
= t_{h}(\ctrq,v\vert_{\Sigma})\; \forall v\in \mbVh(\Omega)$. Besides we have
$\phi(v) = 0\;\forall v\in\Vh(\Omega)$ by hypothesis which rewrites $t_{h}(\ctrq,\ctrp) = 0\;
\forall \ctrp\in \Vh(\Sigma)$ hence $\ctrq\in {\Vh(\Sigma)}^{\perp}$.
This proves \eqref{SystemIntermediaire}.

Reciprocally assume that~\eqref{SystemIntermediaire} holds. According to Lemma~\ref{TraceMatching},
the second equation of~\eqref{SystemIntermediaire} implies that
$u\in \Vh(\Omega)$. Next, taking $v\in \Vh(\Omega)$ in the first equation~\eqref{SystemIntermediaire}
leads to $a(u,v) = \ell(v)$ for all $v\in \Vh(\Omega)$,
which is~\eqref{DiscreteVF}. In conclusion $u$ solves~\eqref{DiscreteVF}. \qed%

\quad\\
The second equation of \eqref{SystemIntermediaire} simply re-writes
$u\vert_{\Sigma}\in\Vh(\Sigma)$. In this saddle point problem, the transmission
conditions stemming from the original problem \eqref{DiscreteVF} are encoded
in the two conditions $u\vert_{\Sigma}\in\Vh(\Sigma)$ (Dirichlet transmission
condition) and $\ctrq\in \Vh(\Sigma)^\perp$ (Neumann transmission condition), and
\eqref{SystemIntermediaire} holds if and only if $(u,\ctrq)\in \mbVh(\Omega)\times
\mbVh(\Sigma)$ satisfies $a(u,v)- t_{h}(\ctrq,v\vert_{\Sigma}) = \ell(v)\;\forall
v\in\mbVh(\Omega)$ and $(u\vert_{\Sigma},\ctrq)\in \Vh(\Sigma)\times {\Vh(\Sigma)}^{\perp}$.

A key  novelty of our analysis consists in using \eqref{CaracTranCond} to reformulate the
two conditions $(u\vert_{\Sigma},\ctrq)\in \Vh(\Sigma)\times {\Vh(\Sigma)}^{\perp}$
equivalently  as $-\ctrq+\iu u\vert_{\Sigma} = \Pi(\ctrq+\iu u\vert_{\Sigma})$.
As a consequence $(u,\ctrq)$ solves~\eqref{SystemIntermediaire}
if and only if it satisfies
\begin{equation}\label{SystemIntermediaire2}
  \begin{aligned}
    & (u,\ctrq)\in \mbVh(\Omega)\times\mbVh(\Sigma)\;\;\text{and}\\
    & a(u,v)-t_{h}(\ctrq,v\vert_{\Sigma}) = \ell(v)\quad \forall v\in\mbVh(\Omega),\\
    & \ctrq -\iu u\vert_{\Sigma} = -\Pi(\ctrq +\iu u\vert_{\Sigma}).
  \end{aligned}
\end{equation}
When it is considered on $\mbVh(\Omega)\times\mbVh(\Omega)$, the sesquilinear form
$a(\cdot,\cdot)$ does not systematically satisfy an inf-sup condition. 
It may admit a non-trivial kernel in certain interior subdomains, which would be an artefact
stemming from domain partitioning only. This is not satisfactory, and this motivates a change
of unknowns $\ctrp = \ctrq-\iu u\vert_{\Sigma}$ which leads to another equivalent formulation
of our problem.

\begin{lem}\quad\\
  Under Assumption \ref{Hypo5}, if  $u\in \Vh(\Omega)$ solves~\eqref{DiscreteVF}
  then there exists a tuple of traces $\ctrp\in\mbVh(\Sigma)$ such that
  \begin{equation}\label{SystemIntermediaire3}
    \begin{aligned}
      & (u,\ctrp)\in \mbVh(\Omega)\times\mbVh(\Sigma)\;\;\text{and}\\
      & \textit{\textcolor{white}{i}(i)}\quad a(u,v)-\iu t_{h}(u\vert_{\Sigma},v\vert_{\Sigma}) = t_{h}(\ctrp,v\vert_{\Sigma})
      + \ell(v)\quad \forall v\in\mbVh(\Omega),\\
      & \textit{(ii)}\quad \ctrp = -\Pi(\ctrp +2\iu u\vert_{\Sigma}).
    \end{aligned}
  \end{equation}
  Reciprocally if $(u,\ctrp)\in \mbVh(\Omega)\times \mbVh(\Sigma)$ solves~\eqref{SystemIntermediaire3},
  then $u\in \Vh(\Omega)$ and it is the unique solution to Problem~\eqref{DiscreteVF}, and 
  $\ctrp+i u\vert_{\Sigma}\in \Vh(\Sigma)^\perp$ and the pair $(u,\ctrq)=(u,\ctrp+i u\vert_{\Sigma})$
  solves \eqref{SystemIntermediaire}.
\end{lem}

\noindent 
Formulation \eqref{SystemIntermediaire3} plays a central role in the
DDM strategy of the present contribution. It contains a volume part (\textit{i})
which expresses wave equations in each subdomain, and a skeleton part (\textit{ii})
that enforces transmission conditions and a (potentially non-local) coupling between
subdomains. Due to the positivity properties of $t_{h}$, the sesquilinear form
$u,v\mapsto a(u,v)-\iu t_{h}(u\vert_{\Sigma},v\vert_{\Sigma})$
satisfies an inf-sup condition.

\begin{lem}[\textbf{Well posedness of local sub-problems}]\label{LocalInfSup}\quad\\
  Under Assumptions \ref{Hypo1}, \ref{Hypo2} and \ref{Hypo5} we have 
  \begin{equation}\label{InfSupCst}
    \beta_{h}:=\mathop{\inf\!\textcolor{white}{p}}_{u\in\mbVh(\Omega)\setminus\{0\}}\,
    \sup_{v\in\mbVh(\Omega)\setminus\{0\}}\frac{\vert\, a(u,v)-\iu
      t_{h}(u\vert_{\Sigma},v\vert_{\Sigma})\,\vert}{\Vert u\Vert_{\mbH^{1}(\Omega)}
      \Vert v\Vert_{\mbH^{1}(\Omega)}}>0. 
  \end{equation}
\end{lem}
\noindent \textbf{Proof:}

Assume that the inf-sup constant above vanishes for some $h>0$.
This implies in particular that there exists $u\in \mbVh(\Omega)\setminus\{0\}$ such that
$a(u,v) -\iu t_{h}(u\vert_{\Sigma},v\vert_{\Sigma}) = 0\;\forall v\in\mbVh(\Omega)$.
Since $\Im m \{\kappa^{2}\}\geq 0$ according to Assumption~\ref{Hypo1}-\textit{(i)},
we have $0 = \int_{\Omega}\mrm{\Im m}\{{\kappa(\bx)}^{2}\}
\vert u\vert^{2} \;\mathrm{d}\bx +\int_{\partial\Omega}\mrm{\Re e}\{\kappa\} \vert u\vert^{2} \;\mathrm{d}\sigma 
+ t_{h}(u\vert_{\Sigma},u\vert_{\Sigma})$ hence $t_{h}(u\vert_{\Sigma},u\vert_{\Sigma}) = 0
\Rightarrow u\vert_{\Sigma} = 0$. From this and Lemma~\ref{TraceMatching}, we conclude that
$u\in \Vh(\Omega)$, and it satisfies $a(u,v) = a(u,v)-\iu t_{h}(u\vert_{\Sigma},v\vert_{\Sigma})
= 0$ for all $v\in\Vh(\Omega)$. This is not possible due to Assumption \ref{Hypo2}. \qed%

\quad\\
Whether or not the inf-sup stability pointed in the previous lemma holds 
uniformly in $h$ is a legitimate question. Because our assumptions on 
$t_{h}(\cdot,\cdot)$ are rather loose though, it is difficult to 
discuss this at present stage. This uniform stability shall be examined 
when we discuss concrete choices of $t_{h}(\cdot,\cdot)$ later on. 

\section{Reduction to a problem on the skeleton}\label{sec:reduction}

Well-posedness of sub-problems allows to introduce local scattering operators
that eliminate volume unknowns, expressing $\ctrp$ in terms of $u$. 
\begin{lem}[\textbf{Scattering operator}]\label{DefOpScattering}\quad\\
  Under Assumptions \ref{Hypo1}, \ref{Hypo2} and \ref{Hypo5}, let
  $\mS_{h}:\mbVh(\Sigma)\to \mbVh(\Sigma)$ refer to the  continuous operator
  defined by $\mS_{h}(\ctrp) = \ctrp + 2\iu w\vert_{\Sigma}$ where
  $w$ is the unique element of  $\mbVh(\Omega)$  satisfying
  $a(w,v)-\iu t_{h}(w\vert_{\Sigma},v\vert_{\Sigma}) = 
  t_{h}(\ctrp,v\vert_{\Sigma})\;\forall v\in\mbVh(\Omega)$.
  Then we have $\Vert \mS_{h}(\ctrp)\Vert_{\Th}\leq \Vert \ctrp\Vert_{\Th}
  \;\forall \ctrp\in \mbVh(\Sigma)$.
\end{lem}
\noindent \textbf{Proof:}

From the definition of $w\in \mbVh(\Omega)$ we deduce that $a(w,w) -
\iu t_{h}(w\vert_{\Sigma},w\vert_{\Sigma}) =  t_{h}(\ctrp,w\vert_{\Sigma})$. The
properties of $a(\cdot,\cdot)$ then guarantee $\Im m\{a(w,w)\}\leq 0$
which leads to the inequality $-\Re e\{\iu t_{h}(\ctrp,w\vert_{\Sigma})\} = 
\Im m\{t_{h}(\ctrp,w\vert_{\Sigma})\}\leq  -\Vert w\vert_{\Sigma}\Vert_{\Th}^{2}$. 
Next, developing the expression of the norm and using the previous inequality,
we obtain  $\Vert \ctrp + 2\iu w\vert_{\Sigma} \Vert_{\Th}^{2} =
\Vert \ctrp\Vert_{\Th}^{2}-4 \Re e\{\iu t_{h}(\ctrp,w\vert_{\Sigma})\}+
4\Vert w\vert_{\Sigma}\Vert_{\Th}^{2}\leq \Vert \ctrp\Vert_{\Th}^{2}$.\qed%

\quad\\
We call $\mS_{h}$ the scattering operator because it implements an ingoing-to-outgoing map
i.e. it takes a tuple of ingoing traces $\ctrp$ as input, solves the associated (discrete) Helmholtz problem
with ingoing Robin traces in each subdomain, and returns $\ctrp + 2\iu w\vert_{\Sigma}$ which is the corresponding
tuple of outgoing Robin traces. The contraction property that we have just established relates to energy
conservation in each subdomain. The scattering operator can be used to eliminate volume unknowns in Formulation
\eqref{SystemIntermediaire3} and re-write it as an equation posed on the skeleton $\Sigma$ only.

\begin{prop}\quad\\
  Under Assumption \ref{Hypo1}, \ref{Hypo2} and \ref{Hypo5}, let $u_{\star}\in \mbVh(\Omega)$
  satisfy $a(u_{\star},v)-\iu t_{h}(u_{\star}\vert_{\Sigma},v\vert_{\Sigma}) = \ell(v)\; \forall v\in\mbVh(\Omega)$
  and set $\ctrf := -2\iu \Pi(u_{\star}\vert_{\Sigma})$. If the pair $(u,\ctrp)\in \mbVh(\Omega)\times\mbVh(\Sigma)$
  solves \eqref{SystemIntermediaire3} then we have
  \begin{equation}\label{OSMFinalForm}
  \begin{aligned}
    & \ctrp\in \mbVh(\Sigma)\;\text{such that}\\
    & (\Id +\Pi\mS_{h})\ctrp = \ctrf.
  \end{aligned}
  \end{equation}
  Reciprocally if $\ctrp\in \mbVh(\Sigma)$ solves \eqref{OSMFinalForm}, and if $u\in \mbVh(\Omega)$
  satisfies $a(u,v)-\iu t_{h}(u\vert_{\Sigma},v\vert_{\Sigma}) = t_{h}(\ctrp,v\vert_{\Sigma})
  + \ell(v)\; \forall v\in\mbVh(\Omega)$, then the pair $(u,\ctrp)\in \mbVh(\Omega)\times\mbVh(\Sigma)$
  solves \eqref{SystemIntermediaire3}.  
\end{prop}

\noindent 
In this proposition, the unique solvability of problems of the form "find $w\in \mbVh(\Omega)$
such that $a(w,v)-it_h(w\vert_{\Sigma},v\vert_{\Sigma}) = $ rhs" is guaranteed by Lemma \ref{LocalInfSup}. 
Equation~\eqref{OSMFinalForm} is a reformulation of~\eqref{SystemIntermediaire3}
as a problem posed on the skeleton $\Sigma$. It is well posed, no matter the right hand side.

\begin{prop}[\textbf{Well posedness of the skeleton formulation}]\label{Isomorphism}\quad\\
  Under Assumptions \ref{Hypo1}, \ref{Hypo2} and \ref{Hypo5}, the operator
  $\Id +\Pi\mS_{h}:\mbVh(\Sigma)\to \mbVh(\Sigma)$ is bijective.
\end{prop}
\noindent \textbf{Proof:}

Since $\dim\,\mbVh(\Sigma)<+\infty$ we only need to prove that $\ker(\Id +\Pi\mS_{h}) = \{0\}$.
Assume $\ctrp\in\mbVh(\Sigma)$ satisfies $(\Id +\Pi\mS_{h})\ctrp = 0$. According
to Lemma~\ref{LocalInfSup}, there exists a unique  $u\in \mbVh(\Omega)$ solving
$a(u,v)-\iu t_h(u\vert_{\Sigma},v\vert_{\Sigma}) =
t_{h}(\ctrp,v\vert_{\Sigma})\;\forall v\in\mbVh(\Sigma)$. Then the pair
$(u,\ctrp)\in \mbVh(\Omega)\times \mbVh(\Sigma)$ solves \eqref{SystemIntermediaire3}
with $\ell \equiv 0$. Since we have established
equivalence between~\eqref{SystemIntermediaire3} and~\eqref{DiscreteVF}, see in
particular Proposition~\ref{EquivalenceFormulation}, we conclude that
$u$ actually belongs to $\Vh(\Omega)$  and solves~\eqref{DiscreteVF}
with right-hand side $\ell \equiv 0$. Hence $u = 0$ according to
Assumption \ref{Hypo2}, and thus $\ctrp = 0$ since
$t_{h}(\ctrp,v\vert_{\Sigma}) = 0\;\forall v\in\mbVh(\Omega)$. \qed%

\quad\\
The  operator $\Id+\Pi\mS_{h}$ admits a special structure ``\textit{identity+contraction}''.
This allows to prove its strong coercivity.

\begin{cor}\label{EstimCoercivite}\quad\\
  Under Assumptions \ref{Hypo1}, \ref{Hypo2} and \ref{Hypo5} we have
  \begin{equation}\label{EstimateEstimCoercivite}
    \begin{aligned}
      & \Re e\{t_{h}(\ctrp,(\Id +\Pi\mS_{h})\ctrp)\}\geq \frac{\gamma_{h}^{2}}{2}\Vert \ctrp\Vert_{\Th}^{2}
      \quad \forall\ctrp\in \mbVh(\Sigma)\\
      & \text{where}\quad \gamma_{h}:= \inf_{\ctrw\in\mbVh(\Sigma)\setminus\{0\}}
    \Vert (\Id + \Pi\mS_{h})\ctrw\Vert_{\Th}/
    \Vert\ctrw\Vert_{\Th}>0.
    \end{aligned}
  \end{equation}
\end{cor}
\noindent \textbf{Proof:}

Due to the contractivity properties given by the Lemmas~\ref{DefOrthoProj}
and~\ref{DefOpScattering}, we have
\begin{equation*}
  \begin{aligned}
    \Vert \ctrp\Vert_{\Th}^{2}\geq \Vert \Pi\mS_{h}(\ctrp)\Vert_{\Th}^{2}
    & = \Vert ((\Id + \Pi\mS_{h})-\Id)\ctrp\Vert_{\Th}^{2}\\
    & = \Vert (\Id + \Pi\mS_{h})\ctrp\Vert_{\Th}^{2} + \Vert \ctrp\Vert_{\Th}^{2}
    - 2\Re e \{ t_{h}(\ctrp,(\Id + \Pi\mS_{h})\ctrp )\}.
  \end{aligned}
\end{equation*}
The result then comes
when eliminating $\Vert \ctrp\Vert_{\Th}^{2}$ from both sides of the
inequality and using the definition of $\gamma_{h}$.  \qed%

\section{Convergent iterative algorithm}\label{sec:algorithm}

Usual iterative methods can be used to compute the solution to Problem~\eqref{SystemIntermediaire3}
or one of the equivalent forms we have obtained for it.
Here we examine the convergence of a Richardson's strategy: 
considering\footnote{The choice of the
  relaxation parameter $r$ follows heuristic considerations.
  If it is chosen to fit explicit calculus for the model geometric configuration
  of two domains and one spherical/circular interface, the value
  $r = 1/\sqrt{2}$ appears to be the optimal choice, see~\cite{MR3989867}.} a fixed relaxation parameter $r\in (0,1)$,
and starting from
given initial data $u^{(0)},\ctrp^{(0)}$, we consider the algorithm
\begin{equation}\label{IterativeAlgo1}
  \begin{aligned}
    & (u^{(n)},\ctrp^{(n)})\in \mbVh(\Omega)\times\mbVh(\Sigma)\;\;\text{and}\\
    & \textit{\textcolor{white}{i}(i)}\quad
    a(u^{(n)},v)-\iu t_{h}(u^{(n)}\vert_{\Sigma},v\vert_{\Sigma}) =
    t_{h}(\ctrp^{(n)},v\vert_{\Sigma})+\ell(v)\quad \forall v\in\mbVh(\Omega),\\
    & \textit{(ii)}\quad \ctrp^{(n)} = (1-r)\ctrp^{(n-1)}- r\Pi(\ctrp^{(n-1)} + 2\iu u^{(n-1)}\vert_{\Sigma}).
  \end{aligned}
\end{equation}
Each step of this algorithm involves two sub-steps. The multi-trace $\ctrp^{(n)}$ should be computed
first through \textit{(ii)} which performs the exchange of information between
subdomains, then $u^{(n)}$ should be computed (possibly in parallel) through \textit{(i)}.

\begin{rmrk}\label{BenefitDiagonality}
Here it appears clearly why choosing diagonal impedance is interesting
(cf Definition~\ref{DiagoImpedance}). Indeed in this case, the sesquilinear
form appearing in \textit{(i)} of~\eqref{IterativeAlgo1} is itself diagonal
$a(u,v)-\iu t_{h}(u\vert_{\Sigma},v\vert_{\Sigma}) = \sum_{j=1}^{\mJ}
a_{\Omega_{j}^{h}}(u,v)-\iu t_{h}^{j}(u\vert_{\Gamma_{j}^{h}}, v\vert_{\Gamma_{j}^{h}})$.
In this situation, solving \textit{(i)} reduces to computing
$u^{(n)} = {(u_{1}^{(n)},\dots,u_{\mJ}^{(n)})}\in \mbVh(\Omega)$ where 
\begin{equation}\label{LocalSubPb}
  \begin{aligned}
      & u_{j}^{(n)}\in \Vh(\Omega_{j}^{h})\quad\text{and}\\
      & a_{\Omega_{j}^{h}}(u_{j}^{(n)},v)-\iu t_{h}^{j}(u_{j}^{(n)}\vert_{\Gamma_{j}^{h}},v\vert_{\Gamma_{j}^{h}}) =
      t_{h}^{j}(p_{j}^{(n)},v\vert_{\Gamma_{j}^{h}})+\ell_{\Omega_{j}^{h}}(v)\quad \forall v\in \Vh(\Omega_{j}^{h})
  \end{aligned}
\end{equation}
where $\ctrp^{(n)} = (p_{1}^{(n)},\dots,p_{\mJ}^{(n)})$. Thanks to the diagonal nature of $t_h(\cdot,\cdot)$, such a problem
must be solved for each $j=1\dots \mJ$ independently. When the impedance is diagonal, solving such problems as
\textit{(i)} of~\eqref{IterativeAlgo1} is then parallel.

This condition on the diagonal nature of the impedance is not required for the present analysis. Yet most of the impedances $t_h(\cdot,\cdot)$
considered in practice agree with the diagonal form of Definition~\ref{DiagoImpedance} and this has important practical
advantages since it is more favorable to distributed-memory parallel implementations where a possible bottleneck lies in
the cost of communication between processors. 
However, certain domain decomposition strategies consider the case
where the impedance operator is not subdomain-wise block diagonal see for example \cite[\S 3.1]{despres:hal-02612368}.
Without assuming that $t_h(\cdot,\cdot)$ is block diagonal though, Algorithm~\eqref{IterativeAlgo1} raises serious computational difficulties.
\end{rmrk}

\quad\\[-10pt] 
The next result shows that the iterative scheme~\eqref{IterativeAlgo1} converges toward
the solution of our initial boundary value problem~\eqref{DiscreteVF} with geometric
convergence.

\begin{thm}[\textbf{Geometric convergence of Richardson algorithm}]\label{ThmConvergenceEstimate}\quad\\
  Under Assumptions \ref{Hypo1}, \ref{Hypo2} and \ref{Hypo5},
  let $\ctrp^{(\infty)}\in\Vh(\Omega)$ refer to the unique solution to~\eqref{OSMFinalForm},
  consider a relaxation parameter $r\in (0,1)$ and define $\gamma_{h}$ as
  in~\eqref{EstimateEstimCoercivite}. Then the iterates $\ctrp^{(n)}$ computed by
  means of~\eqref{IterativeAlgo1} satisfy the estimate
  \begin{equation}\label{ConvergenceEstimate2}
    \frac{\Vert \ctrp^{(n)}-\ctrp^{(\infty)}\Vert_{\Th}}
      {\Vert \ctrp^{(0)}-\ctrp^{(\infty)}\Vert_{\Th}}\leq
      {(1-r(1-r)\gamma_{h}^{2}\,)}^{n/2}.
  \end{equation}
\end{thm}
\noindent \textbf{Proof:}

The $\ctrp^{(n)}$ defined through~\eqref{IterativeAlgo1} satisfy the recurrence
$\ctrp^{(n)} = (1-r)\ctrp^{(n-1)} - r \Pi\mS_{h}(\ctrp^{(n-1)})+r\ctrf$ with
$\ctrf\in\mbVh(\Sigma)$ defined as in~\eqref{OSMFinalForm}. We conclude
that the sequence $\epsilon_{n}:= \ctrp^{(n)} - \ctrp^{(\infty)}$
satisfies $\epsilon_{n} = (1-r)\epsilon_{n-1} - r \Pi\mS_{h}(\epsilon_{n-1})$.
According to Lemma~\ref{DefOrthoProj} and Lemma~\ref{DefOpScattering}, we have
$\Vert \Pi\mS_{h}(\epsilon_{n-1})\Vert_{\Th}\leq \Vert\epsilon_{n-1}\Vert_{\Th}$. From
this we deduce
\begin{equation*}
  \begin{aligned}
    \Vert \epsilon_{n}\Vert_{\Th}^{2}
    & =  \Vert (1-r)\epsilon_{n-1}-r\Pi\mS_{h}(\epsilon_{n-1})\Vert_{\Th}^{2}\\
    & =  (1-r)\Vert \epsilon_{n-1}\Vert_{\Th}^{2} +
    r\Vert \Pi\mS_{h}(\epsilon_{n-1})\Vert_{\Th}^{2}\\
    & \hspace{1.65cm}-r(1-r)\Vert (\Id + \Pi\mS_{h})\epsilon_{n-1}\Vert_{\Th}^{2}\\
    & \leq (1-r(1-r)\gamma_{h}^{2}) \Vert \epsilon_{n-1}\Vert_{\Th}^{2}. 
  \end{aligned}
\end{equation*}
\qed%

\quad\\
Note that the constant \((1-r(1-r)\gamma_{h}^{2})\) is positive, since
\(\gamma_{h}\leq 2\) from the Lemmas~\ref{DefOrthoProj}
and~\ref{DefOpScattering}.
Proposition~\ref{Isomorphism} guarantees that $\gamma_{h}>0$ in the previous result.
This a priori does not discard the possibility that $\lim_{h\to 0}\gamma_{h} = 0$ and 
$\lim_{h\to 0}(1-\gamma_{h}^{2}r(1-r)) = 1$ which would correspond to a deterioration
in the convergence of~\eqref{IterativeAlgo1}. On the other hand, if $\gamma_{h}$ can be proved to
remain bounded away from $0$, this will correspond to \(h\)-uniform geometric convergence.
Through the solution to local sub-problems, the previous theorem
also yields a convergence estimate for $u^{(n)}$.

\begin{cor}\label{VolumicConvergenceRate}\quad\\
  Under Assumptions \ref{Hypo1},\ref{Hypo2} and \ref{Hypo5}, let $(u^{(\infty)},\ctrp^{(\infty)})$ refer to the
  unique solution to~\eqref{SystemIntermediaire3}, and suppose
  that the sequence $(u^{(n)},\ctrp^{(n)})$ has been defined through~\eqref{IterativeAlgo1}.
  Define $\gamma_{h}$ as in~\eqref{EstimateEstimCoercivite}, let $\beta_{h}$ refer to the inf-sup constant
  of Lemma~\ref{LocalInfSup}, and define $\lambda_{h}^{+}$ as in~\eqref{StabilityContinuityModulus}. Then
  we have the estimate
  \begin{equation*}
    \frac{\Vert u^{(n)}-u^{(\infty)}\Vert_{\mbH^{1}(\Omega)}}
      {\hspace{-0.5cm}\Vert\ctrp^{(0)}-\ctrp^{(\infty)}\Vert_{\Th}}\leq
      \frac{\lambda_{h}^{+}}{\beta_{h}}{(1-r(1-r)\gamma_{h}^{2})}^{n/2}.
  \end{equation*}
\end{cor}
\noindent \textbf{Proof:}

Set $e_{n} := u^{(n)}-u^{(\infty)}$ and $\epsilon_{n} := \ctrp^{(n)}-\ctrp^{(\infty)}$.
Combining~\eqref{IterativeAlgo1} with~\eqref{SystemIntermediaire3}
we obtain that $a(e_{n},v)-\iu t_{h}(e_{n}\vert_{\Sigma},v\vert_{\Sigma}) = t_{h}(\epsilon_{n},v\vert_{\Sigma})$
for all $v\in \mbVh(\Omega)$. We have $\Vert v\vert_{\Sigma}\Vert_{\rh}\leq \Vert v\Vert_{\mbH^{1}(\Omega)}$
according to \eqref{HarmonicExtensionMap}-\eqref{MTFNorm}.
Thus, using the inf-sup constant from Lemma~\ref{LocalInfSup}, and the continuity modulus
of $t_{h}$, we obtain $\beta_{h}\Vert e_{n}\Vert_{\mbH^{1}(\Omega)}\leq \lambda_{h}^{+}
\Vert \epsilon_{n}\Vert_{\Th}$. There only remains to apply Theorem~\ref{ThmConvergenceEstimate}. \qed%

\section{Discrete stability}\label{sec:discreteStability}

The inf-sup constant (noted $\gamma_{h}$) of the operator $\Id + \Pi\mS_h$, plays a crucial
role in the bound for the convergence rate provided by Theorem~\ref{ThmConvergenceEstimate}.
In the present section we analyze in more detail  this quantity.

\begin{rmrk}
  Besides $\gamma_{h}$ defined by \eqref{EstimateEstimCoercivite}, the
  forthcoming stability and convergence analysis will repeatedly refer
  to the following constants:
  \begin{itemize}
  \item[$\bullet$] $\alpha_h$ defined through Assumption \ref{Hypo2},
  \item[$\bullet$] $\Vert a \Vert$ defined by \eqref{ContinuityModulusA},
  \item[$\bullet$] $\beta_h$  defined by \eqref{InfSupCst},
  \item[$\bullet$] $\lambda_h^\pm$ defined by \eqref{StabilityContinuityModulus}.
  \end{itemize}
  These constants a priori depend, not only on the mesh width parameter $h$,
  but also on the material characteristics, including the wave number $\kappa$,
  and the geometry of the computational domain.
\end{rmrk}

\noindent 
We first need to introduce ${\mbVh(\Sigma)}^{2}:=\mbVh(\Sigma)\times \mbVh(\Sigma)$ equipped  with the cartesian product
norm defined by
\begin{equation*}
  \Vert (\ctrp,\ctrq)\Vert_{\Th^2}^2
  :=\Vert \ctrp\Vert_{\Th}^{2}+\Vert \ctrq\Vert_{\Th}^{2}.
\end{equation*}
We shall proceed by imitating the analytical approach presented in~\cite{claeys2019new}, which leads to introducing
two subspaces
\begin{equation*}
  \begin{aligned}
      & \mathscr{V}_{h}(\Sigma):=\Vh(\Sigma)\times{\Vh(\Sigma)}^{\perp},\\
    & \mathscr{C}_{h}(\Sigma):=\{\;(\ctru_{\dir},\ctru_{\neu})\in{\mbVh(\Sigma)}^{2}:\; \exists u\in \mbVh(\Omega)\;\text{such that}\;
    u\vert_{\Sigma}=\ctru_{\dir}\\
    & \textcolor{white}{\mathscr{C}_{h}(\Sigma):=\{(\ctru_{\dir},\ctru_{\neu})\in{\mbVh(\Sigma)}^{2},\;\;}
    \text{and}\; a(u,v) = t_{h}(\ctru_{\neu},v\vert_{\Sigma})\;\forall v\in\mbVh(\Omega)\;\}.
  \end{aligned}
\end{equation*}
These are discrete counterparts of the single-trace and Cauchy data spaces considered 
in Section 4 and 6 of~\cite{claeys2019new}. These are two complementary subspaces of the
(discrete) multi-trace space.

\begin{prop}\label{ProjectionBound}\quad\\
  Under Assumptions \ref{Hypo1}, \ref{Hypo2} and \ref{Hypo5}, we have
  $\mbVh(\Sigma)\times \mbVh(\Sigma) = \mathscr{V}_{h}(\Sigma)
  \oplus \mathscr{C}_{h}(\Sigma)$. Moreover, if $\mrm{P}_{h}:\mbVh(\Sigma)\times \mbVh(\Sigma)\to \mathscr{C}_{h}(\Sigma)$
  is the projection onto $\mathscr{C}_{h}(\Sigma)$ with $\ker(\mrm{P}_{h}) = \mathscr{V}_{h}(\Sigma)$, then
  \begin{equation}\label{EstimateProjection}
    \sup_{(\ctrp,\ctrq)\in {\mbVh(\Omega)}^{2}\setminus\{0\}}\frac{\Vert \mrm{P}_{h}(\ctrp,\ctrq)
      \Vert_{\Th^2}}{\Vert(\ctrp,\ctrq)\Vert_{\Th^2}}\leq
      \frac{{(\lambda_{h}^{+})}^{2}+{(2\Vert a\Vert/\lambda_{h}^{-})}^{2}}{\alpha_{h}}.
  \end{equation}
\end{prop}
\noindent \textbf{Proof:}

Assume first that  $(\ctrp,\ctrq)\in\mathscr{V}_{h}(\Sigma)\cap \mathscr{C}_{h}(\Sigma)$.
By definition of \(\mathscr{C}_{h}(\Sigma)\), there exists $u\in \mbVh(\Omega)$ satisfying
$u\vert_{\Sigma}=\ctrp$ and $a(u,v) = t_{h}(\ctrq,v\vert_{\Sigma})\;\forall v\in\mbVh(\Omega)$.
From  $\ctrp\in \Vh(\Sigma)$ and Lemma~\ref{TraceMatching}, we conclude that $u\in\Vh(\Omega)$,
and since $\ctrq\in{\Vh(\Sigma)}^{\perp}$, we have $a(u,v) = 0\forall v\in\Vh(\Omega)$. Hence $u = 0$
according to Assumption \ref{Hypo2}, and thus $\ctrp = \ctrq = 0$. 
This proves that $\mathscr{V}_{h}(\Sigma)\cap \mathscr{C}_{h}(\Sigma) = \{0\}$. 

\quad\\
Next we prove that $\mbVh(\Sigma)\times \mbVh(\Sigma) = \mathscr{V}_{h}(\Sigma)+\mathscr{C}_{h}(\Sigma)$.
Pick $(\ctrp,\ctrq)\in\mbVh(\Sigma)\times \mbVh(\Sigma)$ arbitrarily, and define $\psi$ as
the unique solution to $\psi\in \Vh(\Omega)$ and $a(\psi+\rho_{h}(\ctrp),v) = t_{h}(\ctrq,v\vert_{\Sigma})\ 
\forall v\in\Vh(\Omega)$ according to Assumption \ref{Hypo2}. This function is bounded by 
\begin{equation}\label{Estimate1}
  \Vert \psi\Vert_{\mH^{1}(\Omega)}\leq \big(\frac{\Vert a\Vert}{\alpha_{h}\lambda_{h}^{-}}\big)
  \Vert \ctrp\Vert_{\Th} + \frac{\lambda_{h}^{+}}{\alpha_{h}}\Vert \ctrq\Vert_{\Th}.
\end{equation}
Next rewriting $u := \psi+\rho_{h}(\ctrp)$, we set  $\ctru_{\dir} := u\vert_{\Sigma}$. In addition,
we have $a(u,w) = 0$ for all $w\in\Vh(\Omega)$ satisfying $w\vert_{\Sigma} = 0$ so, according to
Lemma~\ref{LiftingLemma}, we can define $\ctru_{\neu}$ as the unique element of $\mbVh(\Sigma)$
satisfying  $t_{h}(\ctru_{\neu},v\vert_{\Sigma}) = a(u,v)\;\forall v\in \mbVh(\Omega)$ and actually
$t_{h}(\ctru_{\neu},\ctrw) = a(u,\rho_{h}(\ctrw))\;\forall \ctrw\in \mbVh(\Sigma)$. We deduce the estimates
\begin{equation}\label{Estimate2}
  \begin{aligned}
    & \Vert u\Vert_{\mbH^{1}(\Omega)}\leq \Vert \psi\Vert_{\mH^{1}(\Omega)}+
    \Vert\ctrp\Vert_{\Th}/\lambda_{h}^{-},\\
    & \Vert \ctru_{\dir}\Vert_{\Th}\leq \lambda_{h}^{+} \Vert u\Vert_{\mbH^{1}(\Omega)},\\
    & \Vert \ctru_{\neu}\Vert_{\Th}\leq (\Vert a\Vert/\lambda_{h}^{-})\Vert u\Vert_{\mbH^{1}(\Omega)}.
  \end{aligned}
\end{equation}
Now observe that $(\ctru_{\dir},\ctru_{\neu})\in\mathscr{C}_{h}(\Sigma)$ by construction.
Besides we have  $\ctrp - \ctru_{\dir} = \rho_{h}(\ctrp)\vert_{\Sigma} - \ctru_{\dir} =
-\psi\vert_{\Sigma}\in \Vh(\Sigma)$ since $\psi\in\Vh(\Omega)$. Finally we have
$t_{h}(\ctrq-\ctru_{\neu},v\vert_{\Sigma}) = 0$ for all $v\in\Vh(\Omega)$ so that $\ctrq - \ctru_{\neu}\in{\Vh(\Sigma)}^{\perp}$.
We conclude that $(\ctrp,\ctrq) - (\ctru_{\dir},\ctru_{\neu})\in \mathscr{V}_{h}(\Sigma)$ and thus
$(\ctru_{\dir},\ctru_{\neu}) = \mrm{P}_{h}(\ctrp,\ctrq)$. Finally Estimate~\eqref{EstimateProjection} is obtained
by combining~\eqref{Estimate1} with~\eqref{Estimate2} and observing that $\Vert a\Vert/\alpha_{h}\geq 1$ systematically.
\qed%

\quad\\
Next we point that $\mathscr{C}_{h}(\Sigma)$ is closely related to the graph of the
scattering operator $\mS_{h}$, as confirmed by the next lemma.

\begin{lem}\label{CauchySpaceDescription}\quad\\
  Under Assumptions \ref{Hypo1}, \ref{Hypo2} and \ref{Hypo5}, we have\\
  $\mathscr{C}_{h}(\Sigma) = \{(\ctru_{\dir},\ctru_{\neu})\in\mbVh(\Sigma)\times\mbVh(\Sigma):\;
  \ctru_{\neu}+\iu\ctru_{\dir}= \mS_{h}(\ctru_{\neu}-\iu\ctru_{\dir})\}$.
\end{lem}
\noindent \textbf{Proof:}

Take any $(\ctru_{\dir},\ctru_{\neu})\in \mathscr{C}_{h}(\Sigma)$. By definition there exists $u\in \mbVh(\Omega)$
such that $u\vert_{\Sigma} = \ctru_{\dir}$ and $a(u,v) = t_{h}(\ctru_{\neu},v\vert_{\Sigma})$ for all $v\in \mbVh(\Omega)$.
This implies  $a(u,v) -\iu t_{h}(u\vert_{\Sigma},v\vert_{\Sigma})= t_{h}(\ctru_{\neu}-\iu \ctru_{\dir},v\vert_{\Sigma})
\;\forall v\in \mbVh(\Omega)$. Hence we have $\mS_{h}(\ctru_{\neu}-\iu \ctru_{\dir}) = \ctru_{\neu}-\iu \ctru_{\dir}+
2\iu u\vert_{\Sigma} = \ctru_{\neu}+\iu \ctru_{\dir}$ according to the definition of the scattering operator
given by Lemma~\ref{DefOpScattering}.

Reciprocally, assume that $(\ctru_{\dir},\ctru_{\neu})\in \mbVh(\Sigma)$ satisfies $\ctru_{\neu}+\iu\ctru_{\dir}=
\mS_{h}(\ctru_{\neu}-\iu\ctru_{\dir})$. Defining $u$ as the unique element of $\mbVh(\Omega)$ solving 
$a(u,v) - \iu t_{h}(u\vert_{\Sigma}, v\vert_{\Sigma}) = t_{h}(\ctru_{\neu}-\iu\ctru_{\dir}, v\vert_{\Sigma})\;
\forall v\in \mbVh(\Sigma)$, we have $\ctru_{\neu}+\iu\ctru_{\dir}= \ctru_{\neu}-\iu\ctru_{\dir} + 2\iu
u\vert_{\Sigma}\Rightarrow \ctru_{\dir} = u\vert_{\Sigma}$.  From this we also deduce
$a(u,v) = t_{h}(\ctru_{\neu}, v\vert_{\Sigma})$ for all
$v\in\mbVh(\Omega)$. Hence $(\ctru_{\dir},\ctru_{\neu})\in\mathscr{C}_{h}(\Sigma)$.\qed%

\quad\\
The projection we have defined in Proposition~\ref{ProjectionBound} leads to an explicit
expression of the inverse operator ${(\Id+\Pi\mS_{h})}^{-1}$, which we can use to bound the
corresponding inf-sup constant.

\begin{prop}\label{ResolventBound}\quad\\
  Under Assumptions \ref{Hypo1}, \ref{Hypo2} and \ref{Hypo5}, we have
  \begin{equation*}
    \gamma_{h}:=\inf_{\ctrw\in\mbVh(\Sigma)\setminus\{0\}}
    \frac{\Vert (\Id + \Pi\mS_{h})\ctrw\Vert_{\Th}}
      {\Vert\ctrw\Vert_{\Th}}\geq
      \frac{\sqrt{2}\,\alpha_{h}}{{(\lambda_{h}^{+})}^{2}+{(2\Vert a\Vert/\lambda_{h}^{-})}^{2}}.
  \end{equation*}
\end{prop}
\noindent \textbf{Proof:}

Pick an arbitrary $\ctrf\in\mbVh(\Sigma)$ and set $\ctrp_{\dir} = \iu(\Id-\Pi)\ctrf/4$,
$\ctrp_{\neu} = (\Id+\Pi)\ctrf/4$.
Next, considering the projection introduced in Proposition~\ref{ProjectionBound},
define $(\ctru_{\dir},\ctru_{\neu}):= \mrm{P}_{h}(\ctrp_{\dir},\ctrp_{\neu})
\in\mathscr{C}_{h}(\Sigma)$. By construction we have $\ctru_{\dir}-\ctrp_{\dir}\in\Vh(\Sigma)$ and
$\ctru_{\neu}-\ctrp_{\neu}\in{\Vh(\Sigma)}^{\perp}$. Applying Lemma~\ref{DefOrthoProj}, we obtain
\begin{equation*}
  \begin{array}{lrl}
    & -(\ctru_{\neu}-\ctrp_{\neu}) +  \iu(\ctru_{\dir}-\ctrp_{\dir})
    & = \Pi\big( (\ctru_{\neu}-\ctrp_{\neu})+\iu (\ctru_{\dir}-\ctrp_{\dir}) \big)\\
    \iff
    & \ctru_{\neu}-\iu \ctru_{\dir} + \Pi(\ctru_{\neu}+\iu\ctru_{\dir})
    & = (\Id+\Pi)\ctrp_{\neu} -\iu (\Id-\Pi)\ctrp_{\dir}\\
    & & = \lbr (\Id+\Pi)/2\rbr^{2}\ctrf + \lbr (\Id-\Pi)/2\rbr^{2}\ctrf\\
    & & = (\Id+\Pi)\ctrf/2 + (\Id-\Pi)\ctrf/2 = \ctrf.
  \end{array}
\end{equation*}
Applying Lemma~\ref{CauchySpaceDescription}, we conclude that $(\Id +\Pi\mS_{h})(\ctru_{\neu}-\iu \ctru_{\dir})
= \ctrf$ which rewrites $\ctru_{\neu}-\iu \ctru_{\dir} = {(\Id +\Pi\mS_{h})}^{-1}\ctrf$.

Now, observe that, due to the orthogonality of the projectors
$(\Id \pm\Pi)/2$ for the scalar product induced by \(t_{h}\) (Lemma~\ref{DefOrthoProj}), we have
$\Vert (\ctrp_{\dir},\ctrp_{\neu})\Vert_{\Th^2}^{2} = \Vert \ctrp_{\dir}\Vert_{\Th}^{2}
+ \Vert \ctrp_{\neu}\Vert_{\Th}^{2} = \Vert \ctrf\Vert_{\Th}^{2}/4$.
Finally, using Proposition~\ref{ProjectionBound}, we obtain the estimate 
\begin{equation*}
  \begin{aligned}
      \Vert {(\Id +\Pi\mS_{h})}^{-1}\ctrf \Vert_{\Th} 
    & = \Vert \ctru_{\neu}-\iu \ctru_{\dir}\Vert_{\Th} \\
    & \leq \sqrt{2}\Vert (\ctru_{\dir},\ctru_{\neu})\Vert_{\Th^2}
    = \sqrt{2} \Vert \mrm{P}_{h}(\ctrp_{\dir},\ctrp_{\neu})\Vert_{\Th^2}\\
    & \leq \frac{{(\lambda_{h}^{+})}^{2}+ {(2\Vert a\Vert/\lambda_{h}^{-})}^{2}}{\sqrt{2} \alpha_{h}}
    \Vert\ctrf\Vert_{\Th}.
  \end{aligned}
\end{equation*}
\qed%

\begin{rmrk}
Proposition \ref{ResolventBound} is instructive from the perspective
of wave-number dependency. Indeed, there is no hidden constant in the estimate provided for $\gamma_h$.
Since this result holds under Assumptions  \ref{Hypo1}, \ref{Hypo2} and \ref{Hypo5} only,
Proposition \ref{ResolventBound} remains valid when the mesh width $h$ and the wave number $\kappa$
both vary simultaneously. In other words Proposition \ref{ResolventBound} seems useful
for investigating  Optimized Schwarz Methods in the high frequency regime.

In a situation where both  $h$ and $\kappa$ vary, a deterioration of the lower bound for $\gamma_h$
can only come from $\alpha_h,\Vert a\Vert$ or $\lambda_h^\pm$. The $(\kappa,h)$-behavior of $\lambda_h^\pm$
of course totally depends on the choice of impedance $t_h(\cdot,\cdot)$ and nothing can be inferred
without being more specific on this choice. Let us point that the choice of impedance corresponding to
$t_h(\ctrp,\ctrp) = \Vert \ctrp\Vert_{\rh}^2$ from
Example~\ref{ex:schurcomplementimpedance} leads to $\lambda_h^\pm = 1$ which
is clearly $(\kappa,h)$-uniform. The bound $\Vert a\Vert$ can be considered
harmless regarding
$(\kappa,h)$-dependency (see Remark \ref{UniformBoundednessSesquilinearform}). Then only the
$(\kappa,h)$-behavior of $\alpha_h$ remains to be analyzed. Understanding the behavior of this
discrete inf-sup constant for large wave numbers is a non-trivial issue, beyond the
scope of the present article, and this is the subject of active
research at present, see \cite{zbMATH07114221,zbMATH05969646}.
\end{rmrk}

\noindent 
Combining Proposition~\ref{ResolventBound} with Lemma~\ref{EstimCoercivite} or Theorem~\ref{ThmConvergenceEstimate}
obviously yields an estimate for the convergence rate of linear iterative solvers. In particular, when the
impedance is chosen as equivalent to the scalar product associated to $\Vert\cdot\Vert_{\rh}$,
convergence is then uniform with respect to the discretization parameter. 
This motivates to introduce the following condition. We will systematically
state explicitly when this condition is necessary (but it is not assumed to hold in
general).

\begin{cond}[\textbf{\(h\)-uniformly stable impedance}]\label{Hypo6}\quad\\
  The impedance operator is such that\\
  $\lambda_{\star}^{-} :=\liminf_{h\to 0}\lambda_{h}^{-}>0$ and
  $\lambda_{\star}^{+} :=\limsup_{h\to 0}\lambda_{h}^{+}<+\infty$.
\end{cond}

\begin{cor}\label{UniformConvergence}\quad\\
  Under Assumptions \ref{Hypo1}, \ref{Hypo2} and \ref{Hypo5} let $\ctrp^{(\infty)}\in\Vh(\Omega)$
  refer to the unique solution to~\eqref{OSMFinalForm}, consider the relaxation parameter
  $r\in (0,1)$. Assume in addition that material characteristics (in particular $\kappa$)
  are fixed and that Condition~\ref{Hypo6} holds. Then   $0<\liminf_{h\to 0}\gamma_{h}$ and
  for any $0<\gamma_{\star} < \liminf_{h\to 0}\gamma_{h}$ there exists $h_{\star}>0$ such that 
  the iterates $\ctrp^{(n)}$ computed by means of~\eqref{IterativeAlgo1}
  satisfy the estimate
  \begin{equation}\label{ConvergenceEstimate}
    \frac{\Vert \ctrp^{(n)}-\ctrp^{(\infty)}\Vert_{\Th}}
      {\Vert \ctrp^{(0)}-\ctrp^{(\infty)}\Vert_{\Th}}\leq
      {(1-r(1-r)\gamma_{\star}^{2}\,)}^{n/2}\quad \forall h\in (0,h_{\star}),\forall n\geq 0.
  \end{equation}
\end{cor}

\section{Fixed geometric partitions}\label{sec:fixedpartition}

In the present section, to obtain explicit results, we assume that Condition~\ref{Hypo4} holds,
which corresponds to the situation of Figure~\ref{subfig1}. This implies in particular that the number
$\mJ$ of subdomains remains bounded. Besides, we shall also rely on a further condition regarding the impedance
operator \(t_{h}\).

\begin{cond}\label{Hypo7}\quad\\
  There exists a continuous positive definite sesquilinear form $t(\cdot,\cdot):
  \mbH^{1/2}(\Sigma)\times \mbH^{1/2}(\Sigma)\to \CC$ independent of $h$, and two
  constants $c_\pm>0$ independent of $h$, such that $c_-\leq t_{h}(\ctrp,\ctrp)/t(\ctrp,\ctrp)\leq c_+$
  for all $\ctrp\in \mbVh(\Sigma)\setminus\{0\}$ and all $h>0$.
\end{cond}

\noindent 
This condition simply means that the discrete impedance operator is equivalent to a continuous counterpart
that is independent of the mesh. Després impedance (Example \ref{despresimpedance}) and
integral operator based impedances (Example \ref{hypersingularimpedance}) fulfill this condition.
Under the additional condition of $h$-uniform shape regularity \eqref{ShapeRegularity}, the Schur
complement based impedance (Example \ref{ex:schurcomplementimpedance}) also fulfills Condition \ref{Hypo7}.

\begin{lem}\label{UniformInfSup}\quad\\
  Assume that Assumptions~\ref{Hypo1},~\ref{Hypo2} and
  Conditions~\ref{Hypo4},~\ref{Hypo7} are satisfied.
  Then the inf-sup constant $\beta_{h}$ defined in equation~\eqref{InfSupCst} is asymptotically
  $h$-uniformly bounded from below $\beta_{\star} = \liminf_{h\to 0}\beta_{h}>0$.
\end{lem}
\noindent \textbf{Proof:}

We proceed by contradiction, assuming that there exist sequences $h_{n}\to 0$ and
$u_{n}\in \mbV_{h_{n}}(\Omega)$ satisfying $\Vert u_{n}\Vert_{\mbH^{1}(\Omega)} = 1$ and 
\begin{equation}\label{ContradictionHypo}
  \lim_{n\to \infty}\sup_{v\in\mbV_{h_{n}}(\Omega)\setminus\{0\}} \vert a(u_{n},v)-\iu
    t_{h_n}(u_{n},v)\vert/\Vert v\Vert_{\mbH^{1}(\Omega)} = 0.
\end{equation}
Assumptions~\ref{Hypo1} imply that $t_{h_n}(u_n,u_n)\leq \vert \Im m\{ a(u_n,u_n)-\iu t_{h_n}(u_n,u_n)\}\vert$.
Combining Condition \ref{Hypo7} with (\ref{ContradictionHypo}) thus leads to
$t(u_{n},u_{n})\leq t_{h_n}(u_{n},u_{n})/c_-\to 0$. Cauchy-Schwarz inequality applied with $t(\cdot,\cdot)$
then shows 
\begin{equation}\label{ConvergenceImpedanceLimite}
  \lim_{n\to 0}t(u_n,v) = 0\quad \forall v\in \mbH^1(\Omega).
\end{equation}
Next, extracting a subsequence if necessary, we may assume that $u_{n}$ converges toward some $u_{\infty}\in\mbH^{1}(\Omega)$
weakly in $\mbH^{1}(\Omega)$ such that $0=\lim_{n\to \infty}\Vert u_{n}-u_{\infty}\Vert_{\mL^{2}(\Omega)} =
\lim_{n\to \infty}\Vert u_{n}-u_{\infty}\Vert_{\mL^{2}(\partial\Omega)}$. According to \eqref{ConvergenceImpedanceLimite}
we obtain $0 = \lim_{n \to \infty}t(u_n,u_\infty) = t(u_\infty,u_\infty)$. Since $t(\cdot,\cdot)$ is positive definite
on $\mbH^{1/2}(\Sigma)$, we conclude
\begin{equation}\label{ConvergenceImpedanceLimite2}
  t(u_\infty,v) = 0\quad \forall v\in \mbH^1(\Omega).
\end{equation}
Take any $v\in \mbH^{1}(\Omega)$ and let  $v_{n}\in \mbV_{h_{n}}(\Omega)$ refer to its best approximation in
the discrete variational space i.e. $\Vert v-v_{n}\Vert_{\mbH^{1}(\Omega)} = \inf\{\Vert v-w\Vert_{\mbH^{1}(\Omega)},\;
w\in \mbV_{h_{n}}(\Omega)\}$. In particular we have $\Vert v-v_{n}\Vert_{\mbH^{1}(\Omega)}\to 0$. Applying Cauchy-Schwarz
inequality with $t_{h_n}(\cdot,\cdot)$ and using Condition \ref{Hypo7}, we obtain
\begin{equation}\label{ConvergenceImpedanceLimite3}
  \begin{aligned}
    \vert t_{h_n}(u_n,v_n)\vert^2 
    & \leq t_{h_n}(u_n,u_n) t_{h_n}(v_n,v_n)\\
    & \leq c_+ t_{h_n}(u_n,u_n) t(v_n,v_n)\rightarrow 0.
  \end{aligned}
\end{equation}
Weak convergence of $(u_{n})$ together with~\eqref{ContradictionHypo}, \eqref{ConvergenceImpedanceLimite2}
and \eqref{ConvergenceImpedanceLimite3} implies $a(u_{n},v_{n})-\iu t_{h_n}(u_{n},v_{n})\to 0 = a(u_{\infty},v)-\iu t(u_{\infty},v)$.
Since $v$ was chosen arbitrarily in $\mbH^{1}(\Omega)$ we conclude that $u_{\infty}=0$, which implies
$\lim_{n\to \infty}\Vert u_{n}\Vert_{\mL^{2}(\Omega)} = 0$. We finally obtain
$\Vert u_{n}\Vert_{\mbH^{1}(\Omega)}^{2}\leq C\Re e\{a(u_{n},u_{n})\}+
C\Re e\{\int_{\Omega}\kappa^{2}(\bx)\vert u_{n}\vert^{2}\;\mathrm{d}\bx\}+\kappa_{\infty}^2\Vert u_{n}\Vert_{\mL^{2}(\Omega)}^{2}$
where $C = \sup_{\Omega}\vert \mu^{-1}\vert<+\infty$ according to Assumption~\ref{Hypo1}.
Since $\Re e\{a(u_{n},u_{n})\}\to 0$ according to~\eqref{ContradictionHypo} and 
$\Vert u_{n}\Vert_{\mL^{2}(\Omega)} \to 0$, we deduce that $\Vert u_{n}\Vert_{\mbH^{1}(\Omega)}\to 0$
which yields a contradiction. \qed%

\quad\\
From Lemma~\ref{UniformInfSup} and assuming in addition uniform boundedness of
the impedance operator, we easily obtain \(h\)-uniform convergence for the
Richardson algorithm~\eqref{IterativeAlgo1}.

\begin{cor}[\textbf{\(h\)-uniform geometric convergence for Richardson}]\label{UniformVolumicConvergenceRate}\quad\\
  Under Assumptions \ref{Hypo1}, \ref{Hypo2} and \ref{Hypo5}, let $(u^{(\infty)},\ctrp^{(\infty)})$
  refer to the unique solution to~\eqref{SystemIntermediaire3}, and suppose that the sequence
  $(u^{(n)},\ctrp^{(n)})$ has been defined through~\eqref{IterativeAlgo1}. Assume in addition that
  material characteristics (in particular $\kappa$) are fixed; that
  Conditions~\ref{Hypo4},~\ref{Hypo6} and~\ref{Hypo7} are satisfied
  and recall the definition of \(\beta_{\star}\) from Lemma~\ref{UniformInfSup}.
  Then $0<\liminf_{h\to 0}\gamma_{h}$
  and for any $0<\gamma_{\star}<\liminf_{h\to 0}\gamma_{h}$ there exists $h_{\star}>0$ such that 
  \begin{equation}\label{UniformVolumicConvergenceEstimate}
    \frac{\Vert u^{(n)}-u^{(\infty)}\Vert_{\mbH^{1}(\Omega)}}
      {\Vert\ctrp^{(0)}-\ctrp^{(\infty)}\Vert_{\Th}}\leq
      \frac{\lambda_{\star}^{+}}{\beta_{\star}}{(1-r(1-r)\gamma_{\star}^{2})}^{n/2}
      \quad \forall h\in (0,h_{\star}),\forall n\geq 0.
  \end{equation}
\end{cor}

\quad\\
Now let us examine how the estimates of the previous and current sections apply
for the concrete choice of impedance considered in Section~\ref{sec:impedance}.

% \begin{thm}\label{ThmConvergenceEstimate}\quad\\
% \begin{cor}\label{VolumicConvergenceRate}\quad\\
% \begin{cor}\label{UniformConvergence}\quad\\ (Hypo6)
% \begin{lem}\label{UniformInfSup}\quad\\ (Hypo7)
% \begin{cor}\label{UniformVolumicConvergenceRate}\quad\\ (Hypo6 and Hypo7)

\begin{exmpl}[\textbf{Despr{\'e}s impedance}]\label{ConvergenceRateDespresImpedanceFixedPartition}
  In the case of a fixed geometric partition (Condition~\ref{Hypo4}) and fixed
  material characteristics,
  the Despr\'es impedance from Example~\ref{despresimpedance}
  fits within the assumptions (in particular Condition~\ref{Hypo7}) of
  Lemma~\ref{UniformInfSup}.
  However, this choice of impedance violates
  Condition~\ref{Hypo6}, hence the assumptions of
  both results that guarantee \(h\)-uniform convergence rate
  Corollary~\ref{UniformConvergence} and
  Corollary~\ref{UniformVolumicConvergenceRate}.

  It is remarkable that our analysis provides nevertheless explicit upper
  bounds for the convergence factor that give insights on the deterioration of
  the rate of convergence when using this operator.
  Indeed, for this choice of impedance, we have
  $\lambda_{h}^{-} = \mathcal{O}(\sqrt{h})$ and 
  $\limsup_{h\to 0} \lambda_{h}^{+}<+\infty$,
  so that  $\gamma_{h}\geq c h$ for a constant $c>0$ independent of $h$
  according to Proposition~\ref{ResolventBound}.
  Hence, under the conditions of Theorem~\ref{ThmConvergenceEstimate},
  Algorithm~\eqref{IterativeAlgo1} with Despr\'es impedance satisfies the
  following convergence estimate
  \begin{equation}\label{ConvergenceRateDespres}
    \frac{\Vert \ctrp^{(n)}-\ctrp^{(\infty)}\Vert_{\Th}}
    {\Vert \ctrp^{(0)}-\ctrp^{(\infty)}\Vert_{\Th}}\leq 
    {(1-c^2r(1-r) h^2)}^{n/2},
  \end{equation}
  and a similar estimate holds for
  $\Vert u^{(n)}-u^{(\infty)}\Vert_{\mbH^{1}(\Omega)}$,
  see Corollary~\ref{VolumicConvergenceRate}.
  Therefore, although convergence does hold with Despr\'es transmission condition,
  our theory suggests a deterioration in the convergence rate for $h\to 0$.
  The convergence factor upper bound \(\tau:={(1-c^2r(1-r) h^2)}^{1/2}\)
  behaves asymptotically like \(1-\mathcal{O}(h^2)\) as \(h\) goes to \(0\).
  It shall be noticed that convergence factors with such asymptotic behavior
  necessarily induce a quadratic growth of the number of iterations required
  to obtain convergence to a fixed tolerance.
  If we cannot say anything on the sharpness of the upper bound estimate,
  this type of deterioration is somehow visible on our numerical examples
  (which were however not obtained under Condition~\ref{Hypo4}), for instance
  in Figure~\ref{fig:nbit_jacobi_Robin_2D_N4_vsNl_plot}.
\end{exmpl}

\begin{exmpl}[\textbf{Second order differential operator}]\label{ex:comments_sndorder}
  In the case of a fixed geometric partition (Condition~\ref{Hypo4}) and fixed
  material characteristics,
  impedance condition of second order~\eqref{ScndOrder} from
  Example~\ref{sndorderimpedance}
  violates both Conditions~\ref{Hypo6} and~\ref{Hypo7}, hence the assumptions
  of both results that guarantee \(h\)-uniform convergence rate
  Corollary~\ref{UniformConvergence} and
  Corollary~\ref{UniformVolumicConvergenceRate}.

  Again, our analysis provides insights on the deterioration of the rate of
  convergence when using this operator. In this case we have $\liminf_{h\to 0}\lambda_{h}^{-}>0$
  and $\limsup_{h\to 0}\lambda_{h}^{+} = \mathcal{O}(1/\sqrt{h})$ so that, as in the previous
  example, $\gamma_{h}\geq c h$ for a constant $c>0$ independent of $h$, and~\eqref{ConvergenceRateDespres}
  holds. Unlike the preceding example though, we cannot prove $h$-uniform discrete
  local inf-sup stability i.e. $\beta_\star>0$ as given by Lemma~\ref{UniformInfSup}
  hence, in contrast with $\Vert \ctrp^{(n)}-\ctrp^{(\infty)}\Vert_{\Th}$, one cannot
  claim that $\Vert u^{(n)}-u^{(\infty)}\Vert_{\mbH^{1}(\Omega)}$ satisfies an estimate
  similar to \eqref{ConvergenceRateDespres}.  
  Yet, again in this case, the convergence factor upper bound behaves
  asymptotically like \(1-\mathcal{O}(h^2)\) as \(h\) goes to \(0\),
  which would induce a quadratic growth of the number of iterations required
  to obtain convergence to a fixed tolerance if this upper bound was
  describing the actual convergence rate.
  Once again we observe a deterioration in our numerical tests indicating that
  the actual convergence rate behaves similarly to the upper bound with respect
  to the mesh width \(h\), for instance in
  Figure~\ref{fig:nbit_jacobi_Robin_2D_N4_vsNl_plot}.
\end{exmpl}

\begin{exmpl}[\textbf{Integral operator based impedance}]\label{ex:comments_integral}
  In the case of a fixed geometric partition (Condition~\ref{Hypo4}) and fixed
  material characteristics,
  the impedance based on the hypersingular integral
  operator~\eqref{HyperSingular} from Example~\ref{hypersingularimpedance}
  satisfies both Conditions~\ref{Hypo6} and~\ref{Hypo7},
  in contrast with the two previous examples.
  Indeed, in this case, we have
  $0<\liminf_{h\to 0}\lambda_{h}^{-}$ and
  $\limsup_{h\to 0}\lambda_{h}^{+}<+\infty$
  so that $\gamma_{\star} = \liminf_{h\to 0}\gamma_{h}>0$ and we have the
  \(h\)-uniform convergence estimate~\eqref{ConvergenceEstimate} of
  Corollary~\ref{UniformConvergence}.
  Besides, the assumptions of Lemma~\ref{UniformInfSup} are satisfied, which
  yields uniform discrete local inf-sup stability and the
  convergence estimate~\eqref{UniformVolumicConvergenceEstimate} holds for
  $\Vert u^{(n)}-u^{(\infty)}\Vert_{\mbH^{1}(\Omega)}$
  according to Corollary~\ref{UniformVolumicConvergenceRate}.
  Our numerical results do confirm the \(h\)-uniform stability of the
  convergence rate.
\end{exmpl}

\begin{exmpl}[\textbf{Schur complement based impedance}]\label{ex:comments_schur}
  In the case of a fixed geometric partition (Condition~\ref{Hypo4}) and fixed
  material characteristics,
  the Schur complement based impedance~\eqref{SchurComplementBasedImpedance}
  of Example~\ref{ex:schurcomplementimpedance}
  satisfies Condition~\ref{Hypo6}
  since in this case we have $\lambda_h^\pm = 1$.
  Therefore $\gamma_{\star} = \liminf_{h\to 0}\gamma_{h}>0$ and the
  \(h\)-uniform convergence estimate~\eqref{ConvergenceEstimate} of
  Corollary~\ref{UniformConvergence} holds.
  According to Lemma \ref{ScottZhangRobustness}, under the additional condition of $h$-uniform
  shape regularity \eqref{ShapeRegularity}, this operator complies with Condition~\ref{Hypo7} and
  therefore fits the assumptions of Lemma~\ref{UniformInfSup}.  
  This allows to conclude that the convergence estimate~\eqref{UniformVolumicConvergenceEstimate}
  also holds for $\Vert u^{(n)}-u^{(\infty)}\Vert_{\mbH^{1}(\Omega)}$
  according to Corollary~\ref{UniformVolumicConvergenceRate}.
  This is confirmed in our numerical experiments of Section~\ref{sec:numerics}.
  In fact it is remarkable that the best results are systematically associated
  to this choice of impedance.
\end{exmpl}
  
\section{The case of no cross-point}\label{sec:generalization}

The case where the subdomain partition~\eqref{SubdomainPartition} does not
involve any cross-point is an important particular case, so we dedicate the present
section to study this situation. The exchange operator $\Pi$ becomes substantially
simpler in this case. By "absence of cross-point" we mean:
\begin{equation}\label{NoJunction}
  \begin{aligned}
    & \Gamma_{j}^{h}\cap \Gamma_{k}^{h}\cap \Gamma_{p}^{h} = \emptyset \;\;\text{and}\\
    & \Gamma_{j}^{h}\cap \Gamma_{k}^{h}\cap \partial\Omega =\emptyset\\
    & \text{for}\;\; j\neq k,\;k\neq p,\;p\neq j.
  \end{aligned}
\end{equation}
We stress that this "no cross-point" assumption enforces two conditions:
three subdomains cannot be adjacent at any point \textit{and} two subdomains
cannot meet at the physical boundary $\partial\Omega$ of the computational domain.
Examples of such geometric configurations are given in Figure~\ref{Fig2} below.

\begin{figure}[h]
  \centering
  \begin{subfigure}{.64\textwidth}
    \centering
    \includegraphics[height=0.40\textwidth]{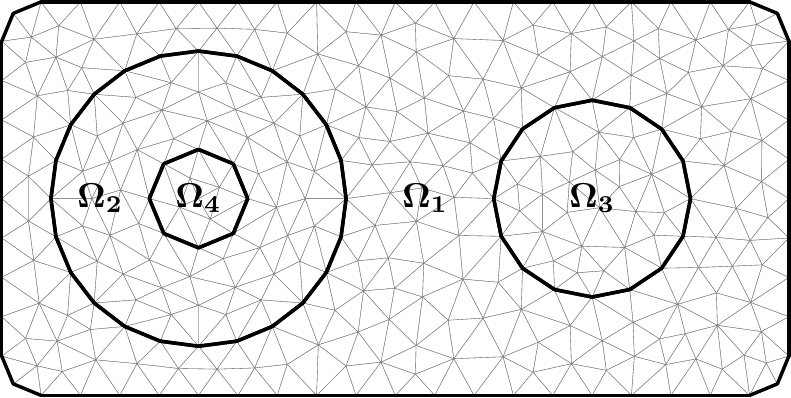}
    \caption{}\label{subfig3}
  \end{subfigure}
  \begin{subfigure}{.32\textwidth}
    \centering
    \includegraphics[height=0.80\textwidth]{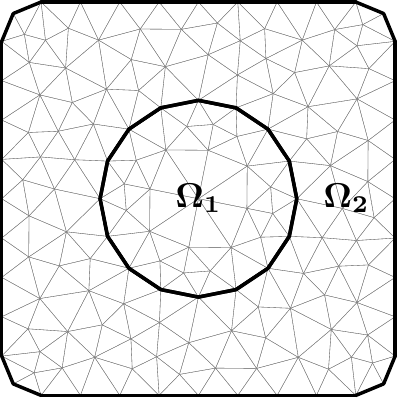}
    \caption{}\label{subfig4}
  \end{subfigure}
  \caption{Examples of partitions without cross-point.}\label{Fig2}
\end{figure}

\noindent 
We wish to study this situation, so we assume all through this
section that~\eqref{NoJunction} holds. Then, each interface is a closed
manifold.  In this case, we can introduce an operator $\mX:\mbVh(\Sigma) \to 
\mbVh(\Sigma)$ consisting in swapping the traces from both sides of each 
interface
\begin{equation}\label{LocalExchange}
  \ctrv = \mX(\ctrw)\;\; \iff\;\;
  \begin{cases}
    v_{j} = w_{k}\;\; \text{on}\;\Gamma_{j}^{h}\cap\Gamma_{k}^{h}\qquad j\neq k,\\
    v_{j} = w_{j}\;\; \,\text{on}\;\Gamma_{j}^{h}\cap\partial\Omega
  \end{cases}
\end{equation}
where $\ctrv = (v_{1},\dots, v_{\mJ})$, $\ctrw = (w_{1},\dots,w_{\mJ})$.
This operator is widely spread in domain decomposition literature.
It is in particular at the core of the Optimized Schwarz Methods (OSM),
see~\cite[Formula (42)]{MR1764190}. The next lemma (simple proof left to
the reader) points a few elementary properties for this operator.

\begin{lem}\quad\\
  If~\eqref{NoJunction} holds then the operator $\mX$ defined
  by~\eqref{LocalExchange} maps continuously $\mbVh(\Sigma)$
  into $\mbVh(\Sigma)$. Besides $\mX^{2} = \Id$ and 
  $(\Id+\mX)/2:\mbVh(\Sigma)\to\Vh(\Sigma)$
  is a projector onto the single-trace space $\Vh(\Sigma)$.
\end{lem}

\noindent 
There are striking similarities shared by both $\Pi$ and $\mX$, see
the definition of $\Pi$ given by Lemma~\ref{DefOrthoProj}. Both operators 
coincide, if and only if $(\Id+\mX)/2$ is orthogonal for the scalar product
$t_{h}(\cdot,\cdot)$, hence if and only if $(\Id+\mX)/2$ is self-adjoint for
the scalar product $t_{h}(\cdot,\cdot)$ which writes equivalently
\begin{equation}\label{SelfAdjointness}
  t_{h}(\mX(\ctrv),\ctrw) =   t_{h}(\ctrv,\mX(\ctrw))
  \quad \forall \ctrv,\ctrw\in \mbVh(\Sigma).
\end{equation}
This condition is satisfied only under certain conditions on the impedance $t_{h}(\cdot,\cdot)$.
The next result provides sufficient conditions for this: roughly speaking, the
impedance operator should be ``symmetric'' with respect to all interfaces
(except the physical boundary).

\begin{prop}\label{prop:xequalspi}\quad\\
    Assume that the impedance \(t_{h}\) is diagonal (cf Definition~\ref{DiagoImpedance})
  and Condition~\eqref{NoJunction} holds.
  Set $\mathfrak{I}:=\{ (j,k)\in {\{1,\dots,\mJ\}}^{2}:\;j<k\;\text{and}\;
  \Gamma_{j}^{h}\cap\Gamma_{k}^{h} \neq\emptyset \}$ and $\Gamma_{0}^{h}:=\partial\Omega$.
  Then $\mX = \Pi$ if there exist scalar products
  $t_{h}^{j,k}:\Vh(\Gamma_{j}^{h}\cap \Gamma_{k}^{h})\times \Vh(\Gamma_{j}^{h}\cap \Gamma_{k}^{h})\to \CC$ such that
  \begin{equation*}
    \begin{array}{l}
        t_{h}(\ctrv,\ctrw) = \sum_{k\in \{1,\ldots,\mJ\}} t_{h}^{0,k}(v_{k},w_{k})
        + \sum_{(j,k)\in \mathfrak{J}} (\,t_{h}^{j,k}(v_{j},w_{j}) + t_{h}^{j,k}(v_{k},w_{k})\,)\\[5pt]
        \forall\; \ctrv = (v_{1},\dots,v_{\mJ}),\;\;\ctrw = (w_{1},\dots, w_{\mJ})\in\mbVh(\Sigma).
    \end{array}
  \end{equation*}
\end{prop}
\noindent  \textbf{Proof:}

From~\eqref{LocalExchange},
for $\ctrv = (v_{1},\dots, v_{\mJ})$ and $\ctrw = (w_{1},\dots,w_{\mJ})$ in \(\mbVh(\Sigma)\),
we obtain the expression
$t_{h}(\mX(\ctrv),\ctrw) = \sum_{k\in \{1,\ldots,\mJ\}}t_{h}^{0,k}(v_{k},w_{k})+
\sum_{(j,k)\in \mathfrak{J}} [t_{h}^{j,k}(v_{j},w_{k}) + t_{h}^{j,k}(v_{k},w_{j})]$.
This proves~\eqref{SelfAdjointness} which, as previously discussed,
shows that $\Pi = \mX$. \qed%

\quad\\
The previous result states that, when there is no cross-point and the impedance
does not couple disjoint interfaces (which is actually a natural choice of impedance),
then $\Pi = \mX$. In this situation, Algorithm~\eqref{IterativeAlgo1} appears to be a
classical Optimized Schwarz Method: this is exactly the algorithm appearing for example
in~\cite[\S 3.3]{MR1291197},~\cite{MR1071633},~\cite[\S 5]{MR1924414}
or~\cite[chap.6]{lecouvez:tel-01444540}.

\quad\\
This shows that~\eqref{IterativeAlgo1} is a true generalization of OSM for the case
where the subdomain partition contains cross-points.

\begin{rmrk}
  From the results of the present section, we deduce in particular that in the case of no cross-point
  and an impedance chosen according to Example \ref{despresimpedance}, our method coincides with 
  the original Despr\'es algorithm introduced in \cite{MR1291197}. The theory from Section
  \ref{sec:algorithm} and \ref{sec:discreteStability} then yields an estimate for the convergence rate
  of this algorithm, see \eqref{ConvergenceRateDespres}. To our knowledge, estimates of the convergence
  rate of Despr\'es algorithm had never been established in such a general setting.
\end{rmrk}

\section{Matrix form of the algorithm}\label{sec:matrixForm}
Coming back to the general situation where the subdomain partition may admit cross-points, 
in this section we will describe in more concrete terms the implementation of the
iterative scheme~\eqref{IterativeAlgo1}, writing all equations in matrix form.
This will help gaining a real insight on the implementation details underlying
the solution strategy we propose. 

\quad\\
First of all, we set a few matrix notations. We assume that the classical shape functions
of $\mathbb{P}_{k}$-Lagrange finite elements are used, and we consider a numbering
of the associated degrees of freedom in each subdomain and on each boundary: for
$\Upsilon = \Omega_{j}^{h}$ or $\Upsilon = \Gamma_{j}^{h}$ for some $j=1,\dots,\mJ$ we define 
\begin{equation*}
  \begin{aligned}
    & N(\Upsilon):=\dim\Vh(\Upsilon)\quad\text{and}\\
    & \Vh(\Upsilon) = \mrm{span}_{k=1\dots N_{\Upsilon}}\{\varphi^{\Upsilon}_{k}\}.
  \end{aligned}
\end{equation*}
Here the $\varphi^{\Upsilon}_{k}$'s refer to the usual $\mathbb{P}_{k}$-Lagrange shape functions associated
to the triangulation in $\Upsilon$. We assume that each shape function on $\Gamma_{j}^{h}$ is obtained by
taking the trace of some shape function on $\Omega_{j}^{h}$.  We also introduce local stiffness matrices
$\bfA_{j}$ with size $N(\Omega_{j}^{h})\times N(\Omega_{j}^{h})$, local impedance matrices $\bfT_{j}$ with
size $N(\Gamma_{j}^{h})\times N(\Gamma_{j}^{h})$, and local trace matrices $\bfB_{j}$ with size
$N(\Gamma_{j}^{h})\times N(\Omega_{j}^{h})$. The entries of these matrices are defined by
\begin{equation*}
  \begin{array}{l}
      {(\bfA_{j})}_{k,l} := a_{\Omega_{j}^{h}}(\varphi_{l}^{\Omega_{j}^{h}},\varphi_{k}^{\Omega_{j}^{h}}),\\[5pt]
      {(\bfT_{j})}_{k,l} := t_{h}^{j}(\varphi_{l}^{\Gamma_{j}^{h}},\varphi_{k}^{\Gamma_{j}^{h}}),
  \end{array}
  \hspace{1cm}
  \left \{
    \begin{array}{l}
        {(\bfB_{j})}_{k,l} := 1\quad \text{if}\;\;\varphi_{k}^{\Gamma_{j}^{h}} = \varphi_{l}^{\Omega_{j}^{h}}\vert_{\Gamma_{j}^{h}},\\[5pt]
        {(\bfB_{j})}_{k,l} := 0\quad \text{otherwise.}
    \end{array}
  \right.
\end{equation*}
Finally we also set
\begin{equation*}
  N(\Sigma) :=\dim\Vh(\Sigma).
\end{equation*}
For each degree of freedom $k = 1,\dots, N(\Sigma)$, set $s(k,j) = 0$
if $k$ does not belong to $\Gamma_{j}^{h}$, and let $s(k,j)$ refer to the number of this degree of freedom local to
$\Gamma_{j}^{h}$ otherwise. The assumptions of conformity that we formulated on the triangulation (see Section~\ref{PbUnderStudy}
and~\ref{GeomPart}) guarantee that we can find a basis of shape functions such that
$\Vh(\Sigma)  = \mrm{span}_{k=1\dots N_{\Sigma}}\{\varphi_{k}^{\Sigma}\}$ where
\begin{equation*}
  \begin{aligned}
      &\varphi_{k}^{\Sigma} = (\varphi_{s(k,1)}^{\Gamma_{1}^{h}},\dots, \varphi_{s(k,\mJ)}^{\Gamma_{\mJ}^{h}})\\
      & \textrm{setting}\;\varphi_{0}^{\Gamma_{j}^{h}}\equiv 0\;\forall j.
  \end{aligned}
\end{equation*}
To keep track of this, we introduce boolean matrices $\bfQ_{j}$ of size $N(\Gamma_{j}^{h})\times N(\Sigma)$
defined by ${(\bfQ_{j})}_{k,l} = 1$ if \(k=s(l,j)\), and ${(\bfQ_{j})}_{k,l} = 0$ otherwise.
These matrices can be used to assemble $\bfT_{\Sigma}$ the Galerkin matrix of the impedance
$t_{h}(\cdot,\cdot)$ restricted to $\Vh(\Sigma)$. It is of size $N(\Sigma)\times N(\Sigma)$
and admits the expression
\begin{equation*}
  \begin{aligned}
      & {(\bfT_{\Sigma})}_{k,l}:=t_{h}(\varphi_{k}^{\Sigma},\varphi_{l}^{\Sigma})\\
    & \bfT_{\Sigma} = \bfQ_{1}^{*}\bfT_{1}\bfQ_{1}+ \cdots +\bfQ_{\mJ}^{*}\bfT_{\mJ}\bfQ_{\mJ}.
  \end{aligned}
\end{equation*}
Finally we also need to introduce local contributions of the right-hand side represented by vectors
$\bfff_{j}$ of size $N(\Omega_{j}^{h})$ defined by $ {(\bfff_{j})}_{k}:= \ell_{\Omega_{j}^{h}}(\varphi_{k}^{\Omega_{j}^{h}})$.
After assembly of the matrices introduced above, and a proper choice of the relaxation parameter $r$
and maximum number of iterations $n_{\max}$, the iterative scheme~\eqref{IterativeAlgo1} takes the
form of Algorithm~\eqref{Algo1} below. The whole algorithm is then parallel except
for the step appearing in Line 10 which ensures coupling between subdomains (see also Remark~\ref{remark1}).

\noindent
\begin{algorithm}[h]
  \caption{}\label{Algo1}
  \begin{algorithmic}[1]
    \For{$j=1,\dots, \mJ$} \Comment{Initialisation}
    \State%
    $\bfp_{j} = 0$ \Comment{size: $N(\Gamma_{j}^{h})$}
    \State%
    $\bfu_{j} = {(\bfA_{j}-\iu \bfB_{j}^{*}\bfT_{j}\bfB_{j})}^{-1}\bfff_{j}$ \Comment{Local solve (size: $N(\Omega_{j}^{h})$)}
    \EndFor%
    \For{$n=1,\dots,n_{\max}$}
    \State%
    $\bfg = 0$ \Comment{size: $N(\Sigma)$}
    \For{$j=1,\dots, \mJ$}
    \State%
    $\bfg = \bfg + \bfQ_{j}^{*}\bfT_{j}(\bfp_{j} + 2\iu\bfB_{j}\bfu_{j})$
    \Comment{Local scattering}
    \EndFor%
    \State%
    $\bfv = \bfT_{\Sigma}^{-1}\bfg$
    \Comment{Global exchange}
    \For{$j=1,\dots, \mJ$}
    \State%
    $\bfp_{j} = \bfp_{j} + 2r (\iu \bfB_{j}\bfu_{j} -\bfQ_{j}\bfv)$
    \State%
    $\bfu_{j} = {(\bfA_{j}-\iu \bfB_{j}^{*}\bfT_{j}\bfB_{j})}^{-1}{(\bfB_{j}^{*}\bfT_{j}\bfp_{j}+\bfff_{j})}$
    \Comment{Local solve (size: $N(\Omega_{j}^{h})$)}
    \EndFor%
    \EndFor%
  \end{algorithmic}
\end{algorithm}

While the theoretical analysis of the Richardson
algorithm~\eqref{IterativeAlgo1} allows to get some deep insight on the
efficiency of the method, such an algorithm is rarely used in practice.
Krylov methods are the preferred choice in real-life applications, in particular
one will typically resort to the \textsc{GMRes} algorithm in our non-symmetric
case.
Importantly, (\(h\)-uniform) geometric convergence of the Richardson algorithm
guarantees (\(h\)-uniform) geometric convergence of its \textsc{GMRes}
counter-part, even the restarted version.

Although other choices are possible, we solve iteratively using \textsc{GMRes}
the linear system given by~\eqref{OSMFinalForm} which features a multi-trace as
unknown. To define the algorithm, it suffices to provide a definition for a
right-hand-side and a matrix-vector product routine.
The right-hand-side is a \(J\)-tuple \((\bfb_{1},\dots,\bfb_{J})\)
and can be computed (offline) according to Algorithm~\ref{Algo2}.
The matrix-vector product procedure, which takes as
input a \(J\)-tuple \((\bfp_{1},\dots,\bfp_{J})\) and
outputs a \(J\)-tuple \((\bfq_{1},\dots,\bfq_{J})\), is given in
Algorithm~\ref{Algo3}.
Notice again here that apart from the computation in Line 7 of
Algorithm~\ref{Algo3} which ensures coupling between subdomains, all
operations are local to the subdomains (see also Remark~\ref{remark1}). 

\noindent
\begin{algorithm}[H]
  \caption{}\label{Algo2}
  \begin{algorithmic}[1]
    \State%
    $\bfg = 0$ \Comment{size: $N(\Sigma)$}
    \For{$j=1,\dots, \mJ$}
    \State%
    $\bfu_{j} = {(\bfA_{j}-\iu \bfB_{j}^{*}\bfT_{j}\bfB_{j})}^{-1}\bfff_{j}$
    \Comment{Local solve (size: $N(\Omega_{j}^{h})$)}
    \State%
    $\bfb_{j} = 2\iu\, \bfB_{j}\bfu_{j}$
    \Comment{size: $N(\Gamma_{j}^{h})$}
    \State%
    $\bfg = \bfg + 2\iu\, \bfQ_{j}^{*}\bfT_{j}\bfB_{j}\bfu_{j}$
    \EndFor%
    \State%
    $\bfv = \bfT_{\Sigma}^{-1}\bfg$
    \Comment{Global exchange}
    \For{$j=1,\dots, \mJ$}
    \State%
    $\bfb_{j} = \bfb_{j} - 2\, \bfQ_{j}\bfv$
    \EndFor%
  \end{algorithmic}
\end{algorithm}

\noindent
\begin{algorithm}[H]
  \caption{}\label{Algo3}
  \begin{algorithmic}[1]
    \State%
    $\bfg = 0$ \Comment{size: $N(\Sigma)$}
    \For{$j=1,\dots, \mJ$}
    \State%
    $\bfu_{j} = {(\bfA_{j}-\iu \bfB_{j}^{*}\bfT_{j}\bfB_{j})}^{-1}{(\bfB_{j}^{*}\bfT_{j}\bfp_{j})}$
    \Comment{Local solve (size: $N(\Omega_{j}^{h})$)}
    \State%
    $\bfq_{j} = - 2\iu\, \bfB_{j}\bfu_{j}$
    \Comment{size: $N(\Gamma_{j}^{h})$}
    \State%
    $\bfg = \bfg + \bfQ_{j}^{*}\bfT_{j}(\bfp_{j}+2\iu\bfB_{j}\bfu_{j})$
    \Comment{Local scattering}
    \EndFor%
    \State%
    $\bfv = \bfT_{\Sigma}^{-1}\bfg$
    \Comment{Global exchange}
    \For{$j=1,\dots, \mJ$}
    \State%
    $\bfq_{j} = \bfq_{j} + 2\, \bfQ_{j}\bfv$
    \EndFor%
  \end{algorithmic}
\end{algorithm}

\section{Numerics}\label{sec:numerics}

In this section, we report on a series of numerical results illustrating the theory of the previous
sections. We emphasize that this section does not aim at evaluating the numerical
performance of our method; this would require fine computational optimization for the
evaluation of the exchange operator $\Pi$, which is beyond the scope of the present contribution.
Discussing the computational cost induced by the exchange operator will be the subject of
another forthcoming article. Our goal here is simply to exhibit numerical confirmation
of our theory.

In all numerical experiments given below, we solve the model
Problem~\eqref{InitialPb} in a domain \(\Omega\) which is either a
disk in 2D or a ball in 3D. Unless stated otherwise, we consider \(\mu\equiv
1\) and the wave number \(\kappa\) is uniform in the domain. The source terms
are taken to be \(f\equiv 0\) and
\(g=(\mu\partial_{\mathbf{n}}-\iu\kappa)u_{\rm inc}\) where \(u_{\rm
inc}(\bx)=\exp(\iu\kappa \mathbf{d}\cdot\bx)\) with
\(\mathbf{d}\) the unit vector in the \(x\) direction.

We provide numerical results obtained for the Richardson
algorithm~\eqref{IterativeAlgo1} as well as results obtained with a restarted
\textsc{GMRes} algorithm. In all our numerical experiments, the relaxation
parameter of the Richardson algorithm is \(r=0.5\) and 
\textsc{GMRes} is restarted every \(20\) iterations.
We provide various tables reporting the number of iterations required to
achieve a tolerance of \(10^{-8}\) for the relative error defined at the
iteration \(n\) as
\begin{equation}\label{eq:relative_error}
  (\text{relative error})^{2} = 
  \frac{\sum_{j=1\dots \mJ}\|u^{(n)}-u^{(\infty)}\|_{\mH^1(\Omega_{j}^{h})}^{2}}
  {\sum_{j=1\dots \mJ}\|u^{(0)}-u^{(\infty)}\|_{\mH^1(\Omega_{j}^{h})}^{2}},
\end{equation}
where \(u^{(n)}\) is the volume solution at iteration \(n\), 
\(u^{(0)}\) is the initial volume solution (taken to be zero in practice) and 
\(u^{(\infty)}\) is the exact discrete volume solution of the full
(undecomposed) problem.
The choice of this volume (energy) norm has the important benefit of being
independent of the choice of impedance or mesh partition. Finally, we
stress that the criterion for reaching convergence does not rely on the
residual of the linear system that is solved.
In all test runs that were performed, the convergence was stopped if a maximum
number of \(10^5\) iterations was not enough to reach the set tolerance.

The impedance operators tested are:
the Despr\'es impedance operator of Example~\ref{despresimpedance}, denoted by
\({\rm M}\) with parameter \(\kappa_R=\kappa\);
the second order impedance operator of Example~\ref{sndorderimpedance}, denoted
by \({\rm K}\) with parameters \(a=\frac{1}{2\kappa}\) and \(b=\kappa\);
the (hypersingular) boundary integral operator given in
Example~\ref{hypersingularimpedance}, denoted by
\({\rm W}\) with parameters \(a=\kappa^2\) and \(\delta=\frac{1}{\kappa}\);
and the Schur complement based operator of
Example~\ref{ex:schurcomplementimpedance}, denoted by \(\Lambda\).

\begin{figure}[h]
    \centering
    \begin{subfigure}{.3\textwidth}
        \includegraphics[width=\textwidth]{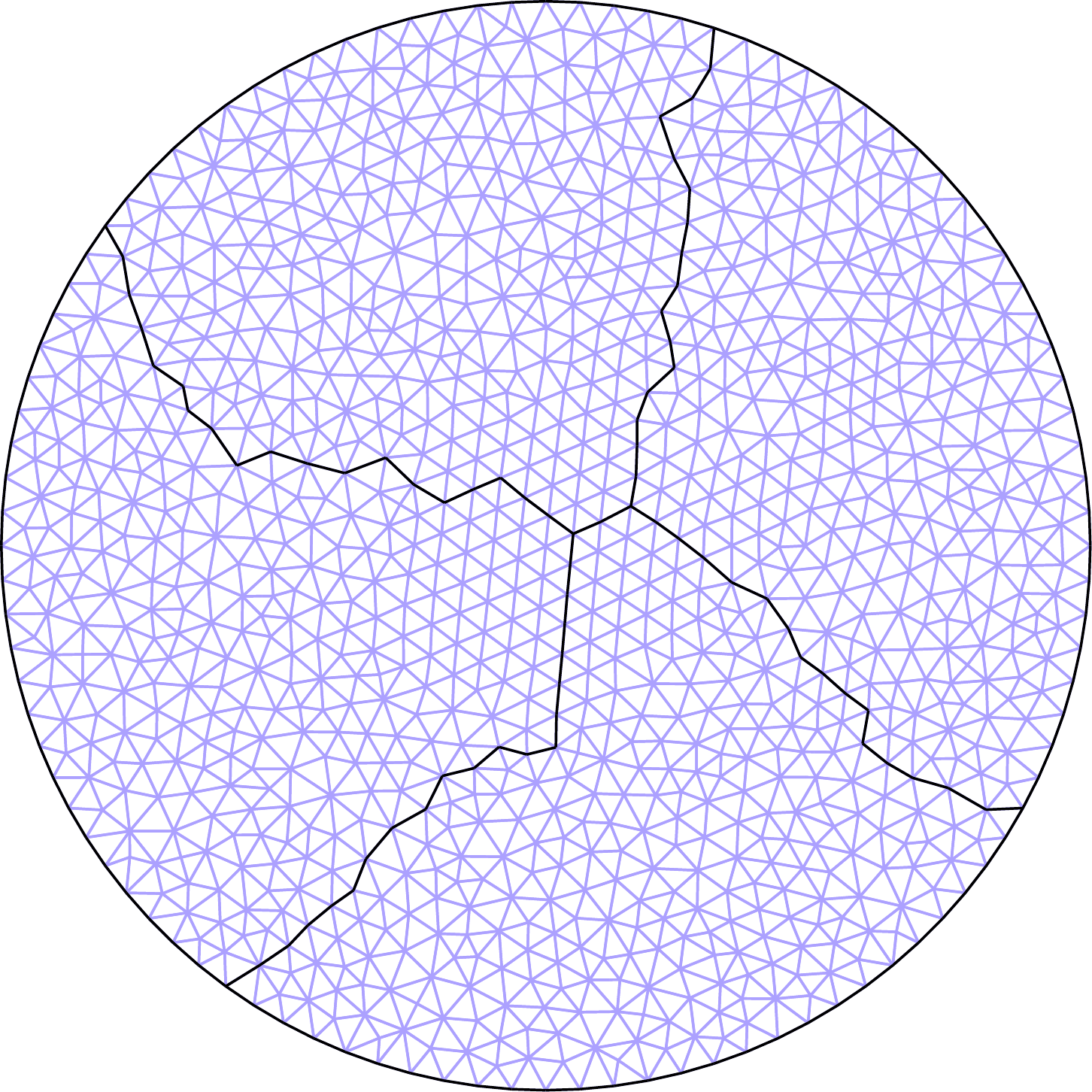}%
        \caption{\(\mJ=4\).}\label{fig:xpts_N4}
    \end{subfigure}
    \begin{subfigure}{.3\textwidth}
        \includegraphics[width=\textwidth]{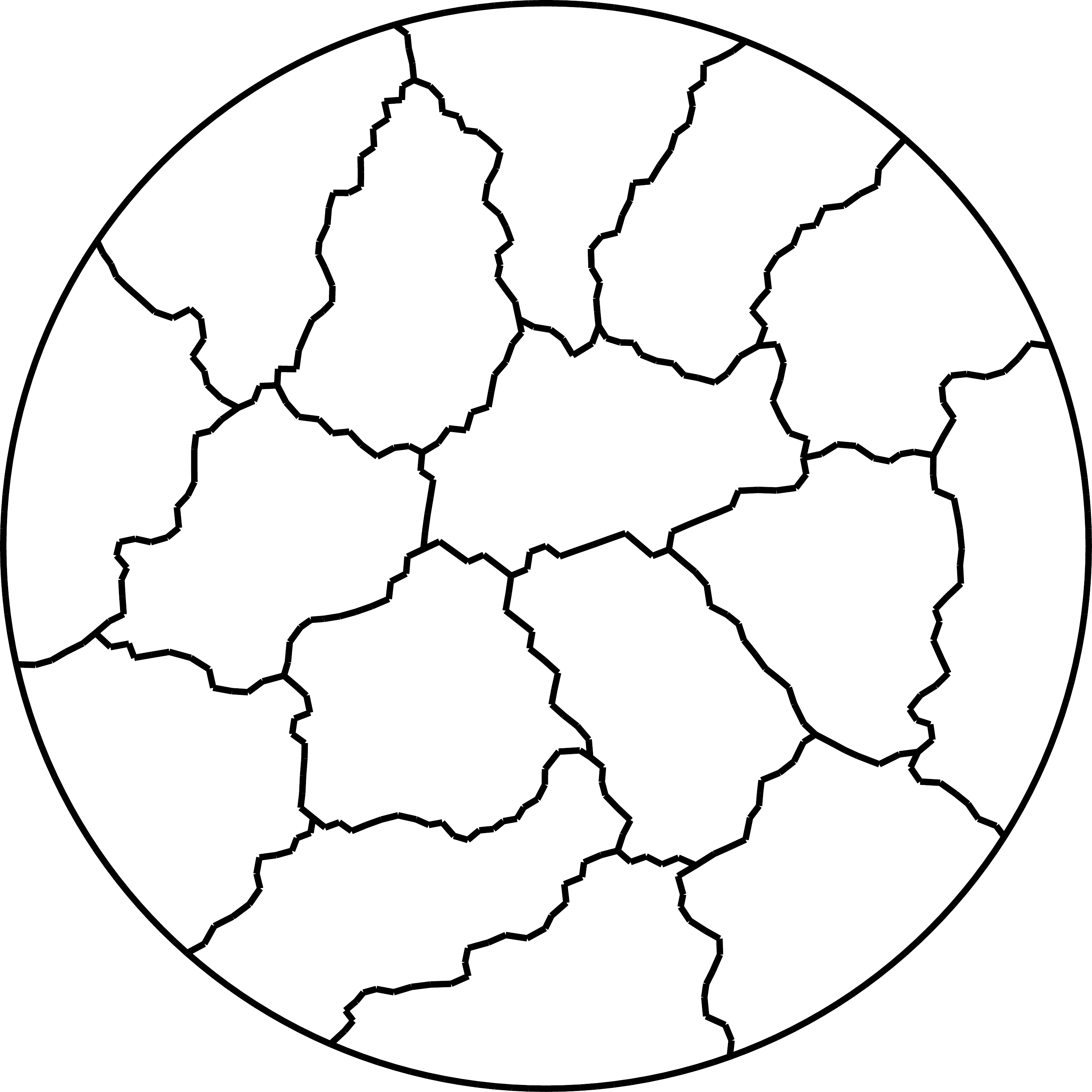}%
        \caption{\(\mJ=16\).}\label{fig:xpts_N16}
    \end{subfigure}
    \begin{subfigure}{.3\textwidth}
        \includegraphics[width=\textwidth]{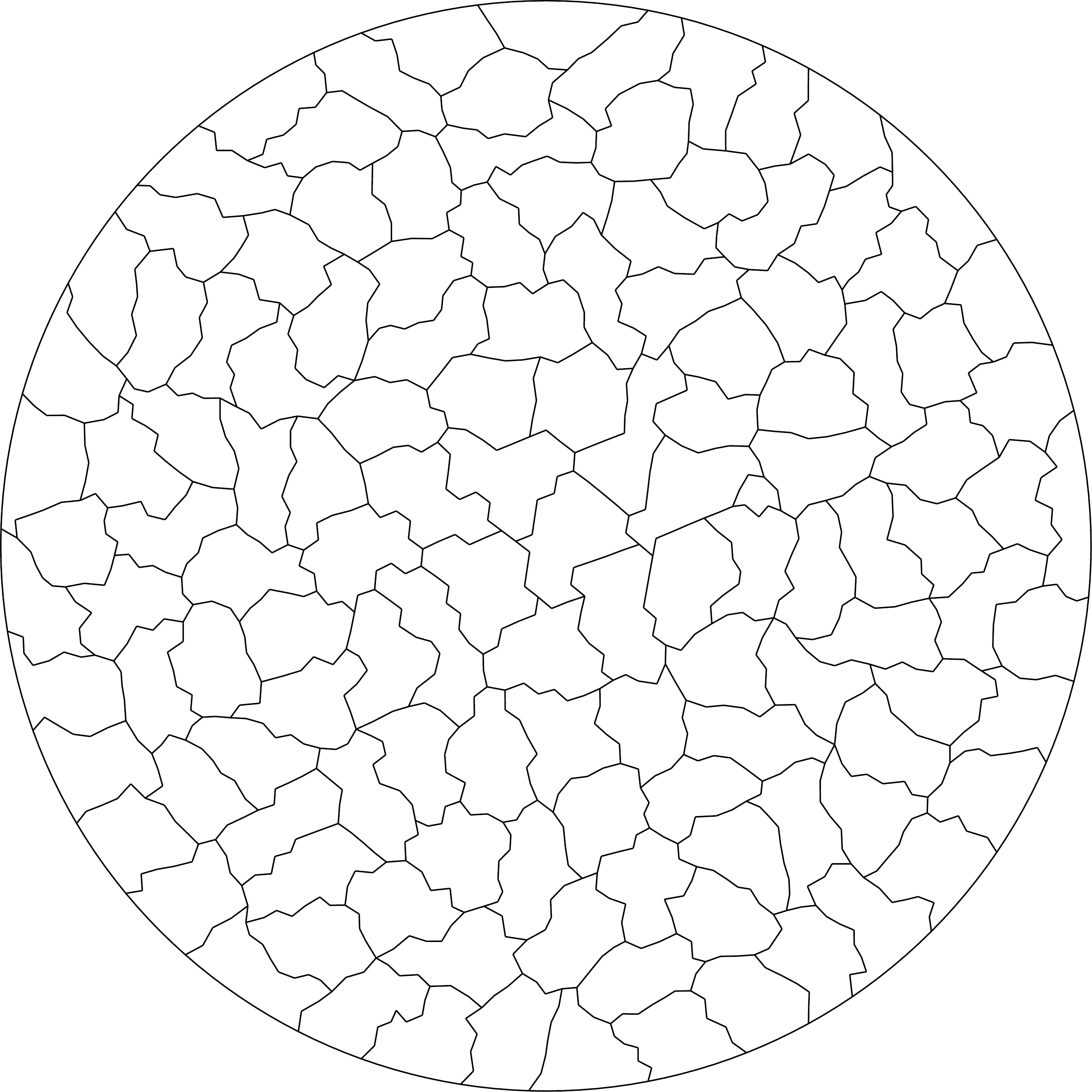}%
        \caption{\(\mJ=128\).}\label{fig:xpts_N128}
    \end{subfigure}\caption{Examples of mesh partitions.}\label{fig:xpts_Nomega}
\end{figure}

The research code that was used to run the tests was developed specifically to
test the method and uses \(\mathbb{P}_{1}\)-Lagrange finite elements. It is
written in \textsc{Julia}~\cite{Bezanson2017} and was validated on standard
scattering test cases. The meshes are generated by
\textsc{Gmsh}~\cite{Geuzaine2009} and partitioned using
\textsc{Metis}~\cite{Karypis1999} through the \textsc{Julia} API\@.
The integral operator matrices are computed thanks to the \textsc{BemTool}
library\footnote{\href{https://github.com/xclaeys/BemTool}{https://github.com/xclaeys/BemTool}},
written in \textsc{C++}.
% The tests were performed on a 8-cores \textit{Intel XEON W-2145} at \(3.7\) GHz
% equipped with 256 Go of RAM\@.

\subsection{Influence of typical mesh size}

We present a first test case consisting of a disk of radius \(R=1\) split
roughly (using a mesh partitioner) in four quarters, see
Figure~\ref{fig:xpts_N4}. The interest of this test case is the presence of
pure interior cross-points where three domains share a common vertex.

\begin{figure}[h]
    \centering
    \begin{subfigure}{.49\textwidth}
        \centering
        \includegraphics[width=\textwidth]{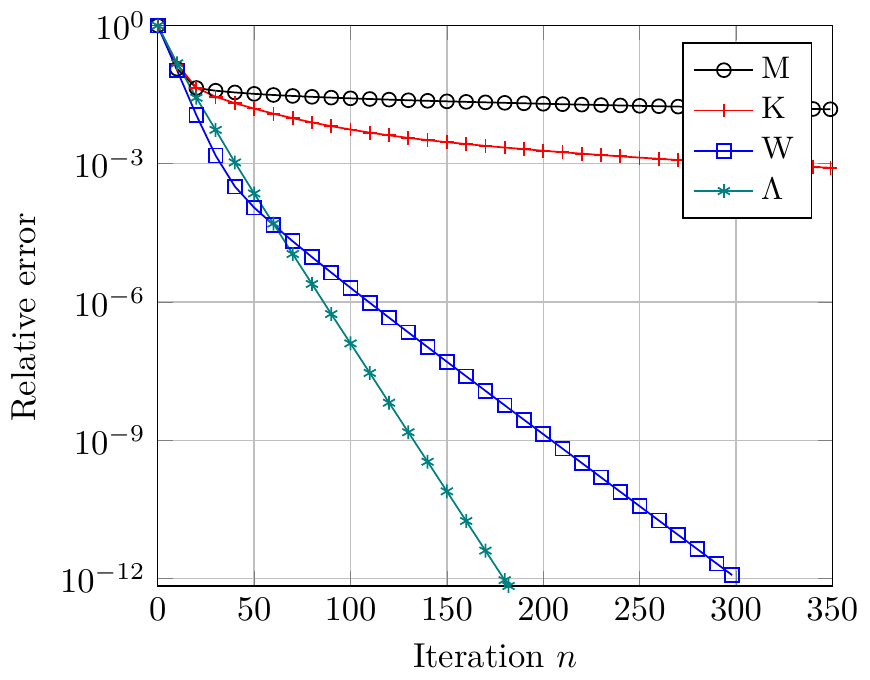}%
        \caption{Richardson algorithm}\label{fig:xpts_cvplot_jacobi_Robin}
    \end{subfigure}
    \begin{subfigure}{.49\textwidth}
        \centering
        \includegraphics[width=\textwidth]{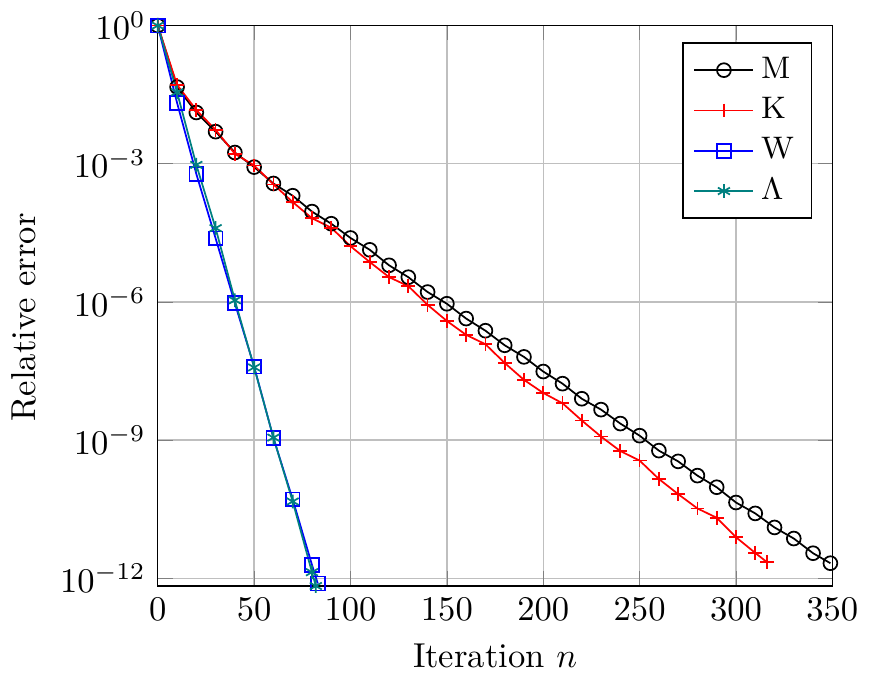}%
        \caption{\textsc{GMRes} algorithm}\label{fig:xpts_cvplot_gmres_Robin}
    \end{subfigure}\caption{An example of convergence history. Fixed parameters
    \(\kappa=5\), \(\mJ=4\), \(N_{\lambda}=2\pi/(\kappa h)=40\),
    2D, disk of radius \(R=1\).}\label{fig:xpts_cvplot_Robin}
\end{figure}

The full convergence history of the relative \(\mathbb{H}^1\)
error~\eqref{eq:relative_error} for the Richardson and \textsc{GMRes}
algorithms are provided for this test case in
Figure~\ref{fig:xpts_cvplot_Robin} as an illustrative example of typical
convergence.

\begin{figure}[h]
    \centering
    \begin{subfigure}{.49\textwidth}
        \centering
        \includegraphics[width=\textwidth]{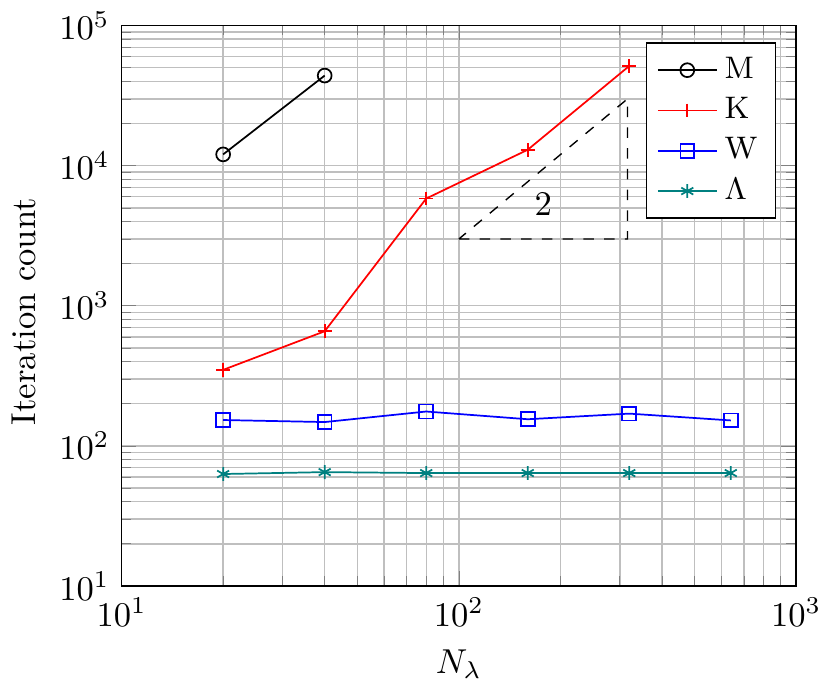}%
        \caption{Richardson algorithm}\label{fig:nbit_jacobi_Robin_2D_N4_vsNl_plot}
    \end{subfigure}
    \begin{subfigure}{.49\textwidth}
        \centering
        \includegraphics[width=\textwidth]{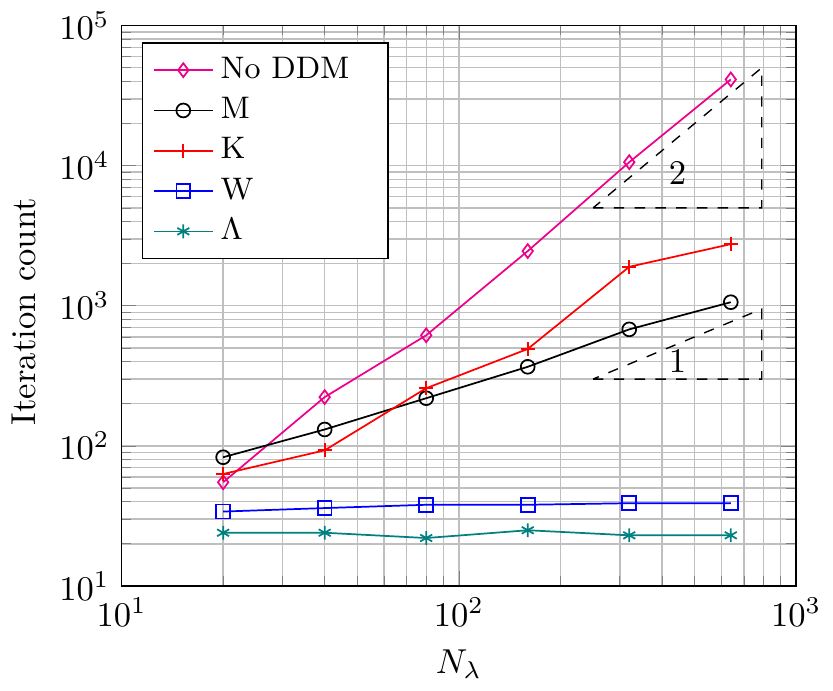}%
        \caption{\textsc{GMRes} algorithm}\label{fig:nbit_gmres_Robin_2D_N4_vsNl_plot}
    \end{subfigure}\caption{Number of iterations with respect to the number of
    mesh points per wavelength \(N_{\lambda}=2\pi/(\kappa h)\). Fixed parameters \(\kappa=1\),
    \(\mJ=4\), 2D, disk of radius \(R=1\).}\label{fig:nbit_Robin_2D_N4_vsNl_plot}
\end{figure}

We report the number of iterations to reach
convergence with respect to mesh refinement in
Figure~\ref{fig:nbit_Robin_2D_N4_vsNl_plot}
for the Richardson and \textsc{GMRes} algorithms.
The refinement of the mesh is indicated by the
number of points per wavelength \(N_{\lambda}=2\pi/(\kappa h)\) which is inversely proportional
to the typical mesh size \(h\).
In Figure~\ref{fig:nbit_gmres_Robin_2D_N4_vsNl_plot} we
also report the number of GMRes iterations that are required to achieve the
same error to solve the full (undecomposed) linear system (line plot labeled `No
DDM'). We see that this iteration count has a growth which is approximately
quadratic with respect to \(N_{\lambda}\), illustrating the deterioration of
the matrix conditioning as the mesh is refined.

For the local operators \({\rm M}\) and \({\rm K}\) the
convergence is not uniform with respect to mesh refinement and a large number
of iterations is required to reach the set tolerance.
The growth of the iteration count appears to be quasi quadratic with respect to
\(N_{\lambda}\) for the Richardson algorithm and quasi-linear for the
\textsc{GMRes} algorithm.
This agrees with the comments on the analytical estimates from
Example~\ref{ConvergenceRateDespresImpedanceFixedPartition}
and Example~\ref{ex:comments_sndorder} which therefore seem sharp.
For small mesh size the convergence may not even be reached within \(10^5\)
iterations.
In contrast, the non-local operators \({\rm W}\) and \(\Lambda\) exhibit
uniform convergence in all cases, with a very moderate number of iterations
required to reach the set tolerance.
This is also in agreement with the comments on the analytical estimates from
Example~\ref{ex:comments_integral} and
Example~\ref{ex:comments_schur}.

\begin{figure}[h]
    \centering
    \centering
    \begin{subfigure}{.49\textwidth}
        \centering
        \includegraphics[width=\textwidth]{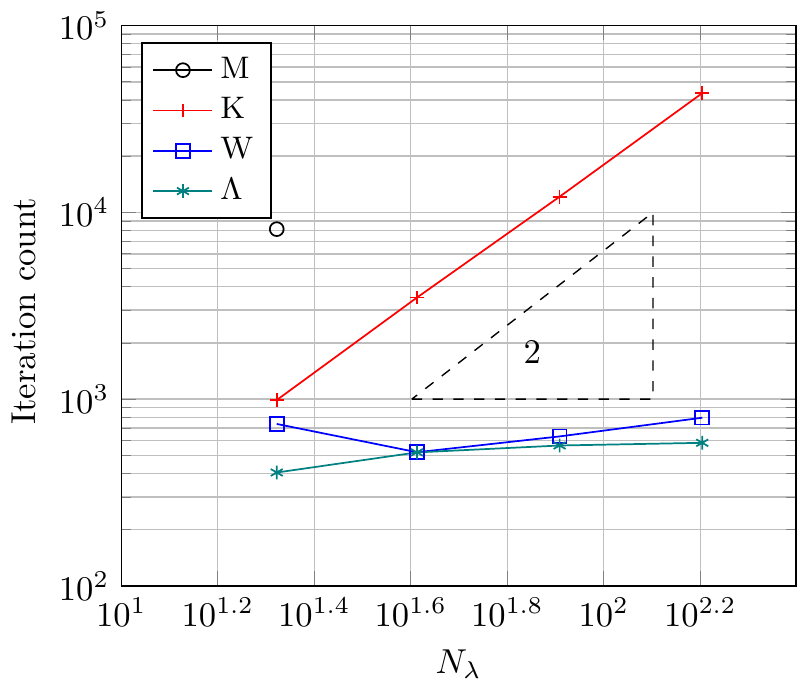}%
        \caption{Richardson algorithm}\label{fig:nbit_jacobi_Robin_3D_N8_vsNl_plot}
    \end{subfigure}
    \begin{subfigure}{.49\textwidth}
        \centering
        \includegraphics[width=\textwidth]{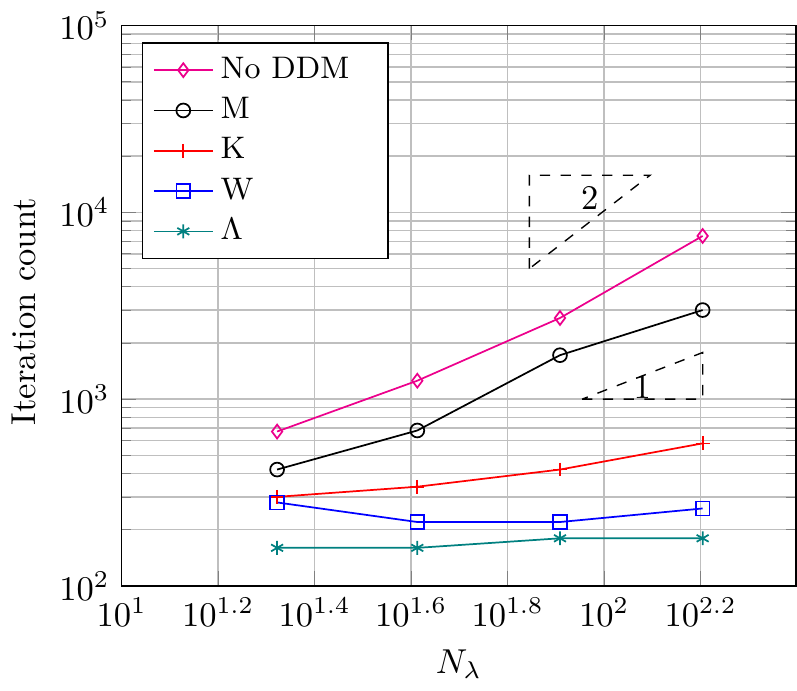}%
        \caption{\textsc{GMRes} algorithm}\label{fig:nbit_gmres_Robin_3D_N8_vsNl_plot}
    \end{subfigure}\caption{Number of iterations with respect to the number of
    mesh points per wavelength \(N_{\lambda}=2\pi/(\kappa h)\). Fixed parameters \(\kappa=1\),
    \(\mJ=8\), 3D, sphere of radius \(R=1\).}\label{fig:nbit_Robin_3D_N8_vsNl_plot}
\end{figure}

We also provide some numerical results obtained in 3D. The domain \(\Omega\) is
now a ball of radius \(R=1\) partitioned into eight subdomains, which generate
interior cross-point curves where three domains share common edges.
Figure~\ref{fig:nbit_Robin_3D_N8_vsNl_plot} reports the iteration count with
respect to mesh refinement. Again in this case, we clearly identify the
non-uniformity of the convergence for the local operators \({\rm M}\) and
\({\rm K}\) while the non-local operators \({\rm W}\) and \({\rm \Lambda}\)
exhibit \(h\)-uniform convergence.

\subsection{Influence of the wave number}

\begin{figure}[h]
    \centering
    \begin{subfigure}{.49\textwidth}
        \centering
        \includegraphics[width=\textwidth]{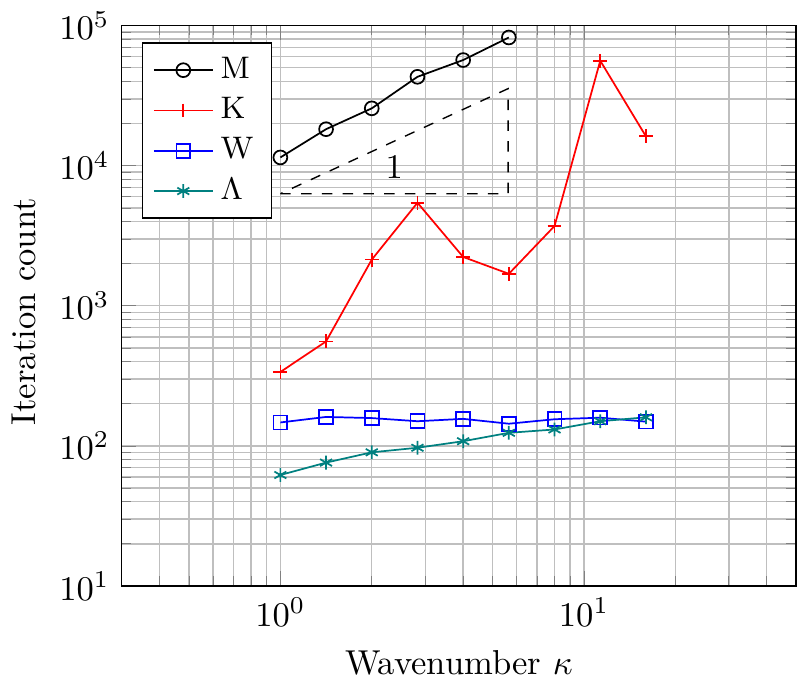}%
        \caption{Richardson algorithm}\label{fig:nbit_jacobi_Robin_2D_N4_vsk_plot}
    \end{subfigure}
    \begin{subfigure}{.49\textwidth}
        \centering
        \includegraphics[width=\textwidth]{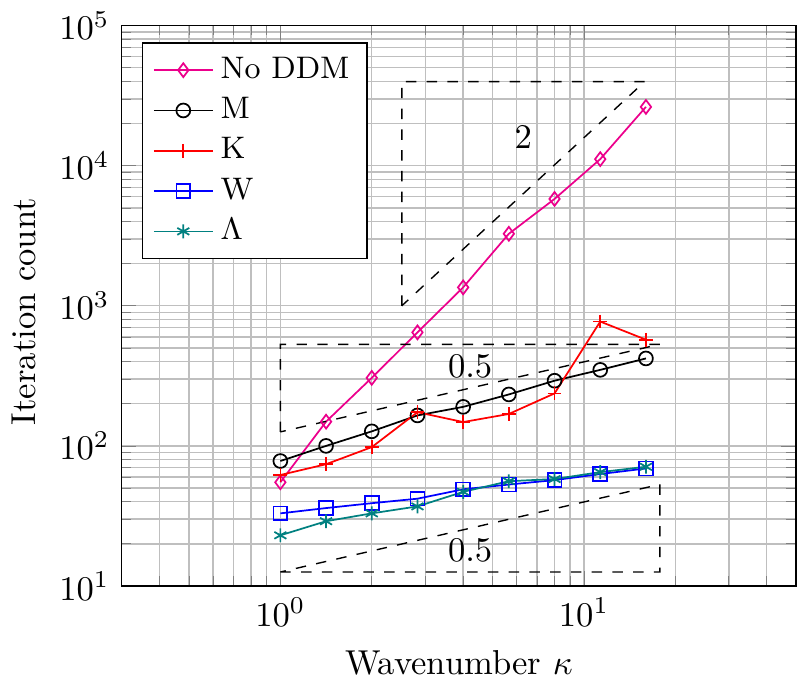}%
        \caption{\textsc{GMRes} algorithm}\label{fig:nbit_gmres_Robin_2D_N4_vsk_plot}
    \end{subfigure}\caption{Number of iterations with respect to the wave number
    \(\kappa\). Fixed parameters \(\mJ=4\), \(h^2\kappa^3 = (2\pi/20)^2\),
    2D, disk of radius \(R=1\).}\label{fig:nbit_Robin_2D_N4_vsk_plot}
\end{figure}

For the two-dimensional case, we now report the dependency of the iteration
count with respect to the wave number \(\kappa\), see
Figure~\ref{fig:nbit_Robin_2D_N4_vsk_plot}. To avoid pollution induced by
phase error, the mesh width varies so as to maintain the relation
$h^2\kappa^3 = (2\pi/20)^2$ i.e. $N_{\lambda}=20$ points per wavelength at wave number
$\kappa = 1$.

As the wave number \(\kappa\) increases, the discrete (as well as the
continuous) problem becomes more difficult to solve. This is indicated again by
the increase in the iteration count of the GMRes algorithm for the undecomposed
problem (line plot labeled `No DDM' in Figure~\ref{fig:nbit_gmres_Robin_2D_N4_vsk_plot}).
On the other hand, for all the impedance operators under study, we notice
a sub-linear growth of the number of iteration with respect to \(\kappa\).
The growth of iteration count is especially moderate for non-local impedances.
Surprisingly, with $\mathrm{W}$ as impedance, the iteration count seems stable
for Richardson's algorithm but not for GMRes. We have no explanation for this fact.

\subsection{Influence of the number of subdomains}

We finally study the dependency of the method with respect to
the number of subdomains \(\mJ\) of the mesh partition.

\begin{figure}[h]
    \centering
    \includegraphics[width=0.49\textwidth]{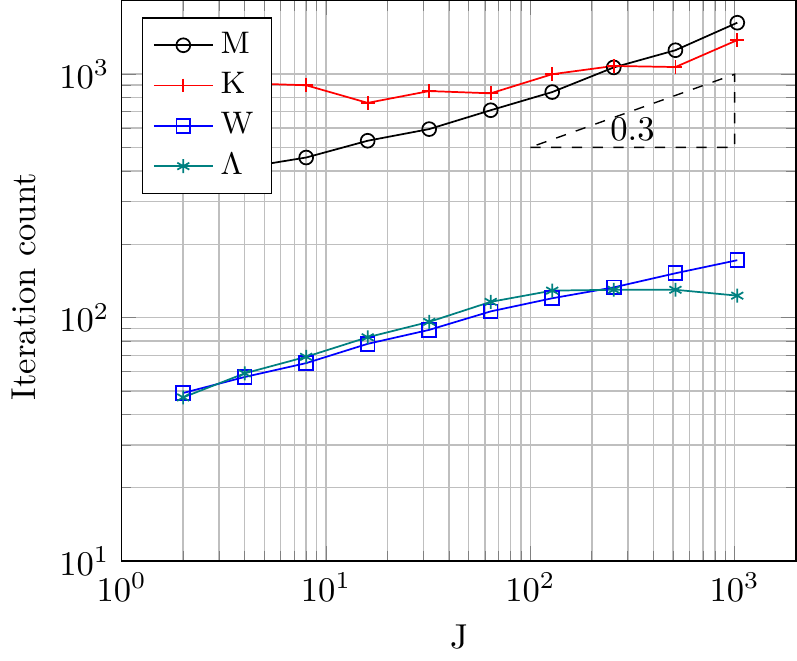}%
    \caption{Number of iterations with respect to the number of subdomains
    \(\mJ\) for a fixed size problem. Fixed parameters \(\kappa=2\), \(N_{\lambda}=100\),
    2D, disk of radius \(R=4\). \textsc{GMRes} algorithm.}\label{fig:nbit_Robin_2D_vsn_plot}
\end{figure}

In the first set of results, we study the influence of the number of subdomains \(\mJ\)
on the iteration count for a problem with fixed size.
Figure~\ref{fig:nbit_Robin_2D_vsn_plot} reports the iteration count with
respect to \(\mJ\) varying from \(2\) to \(1024\) subdomains for a 2D
configuration only.
One can notice a sub-linear increase in the number of iterations to get to a
converged solution. Notice that in this case the undecomposed linear system is
kept the same. Hence, the fact that the discrete problem gets harder is a pure
artificial effect of the DDM.\@ Interestingly, we see that the number of
iterations levels out for the coercive DtN operator, in a regime where the size
of the sub-problems gets really small compared to the wavelength of the
problem.

\begin{figure}[h]
    \centering
    \begin{subfigure}{.49\textwidth}
        \centering
        \includegraphics[height=0.85\textwidth]{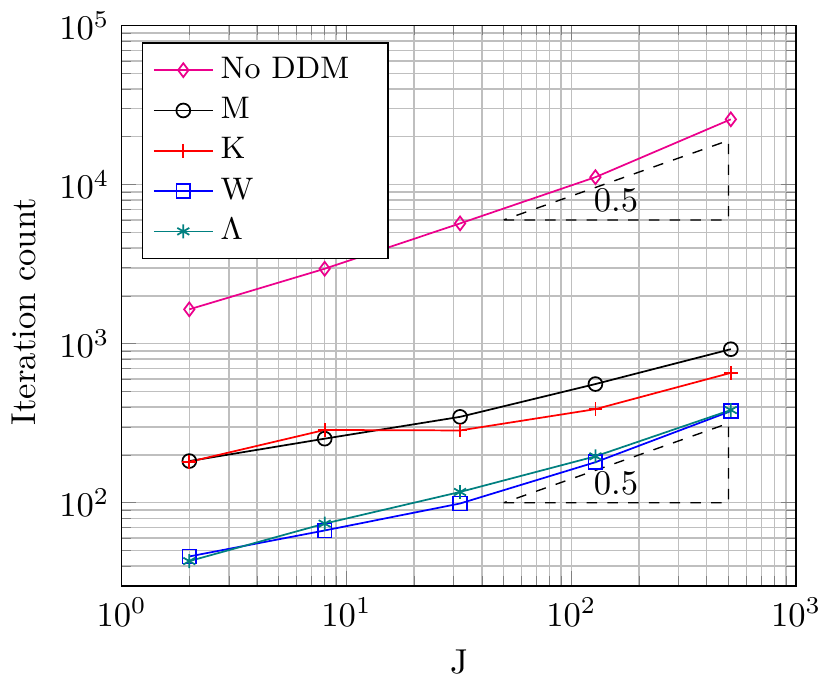}%
        \caption{Fixed parameters \(\kappa=5\), \(N_{\lambda}=40\), 2D}\label{fig:nbit_gmres_Robin_2D_vsnweak_plot}
    \end{subfigure}
    \begin{subfigure}{.49\textwidth}
        \centering
        \includegraphics[height=0.85\textwidth]{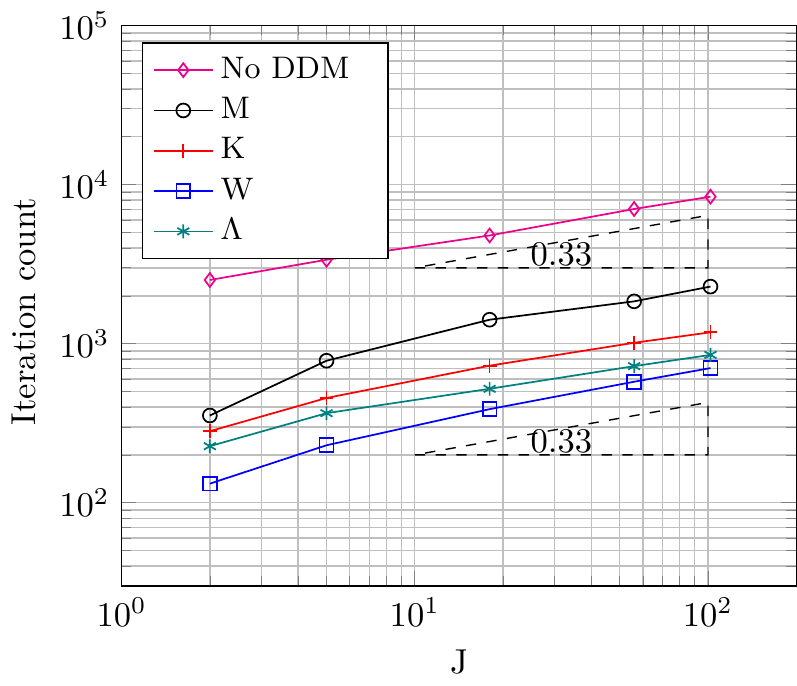}%
        \caption{Fixed parameters \(\kappa=2\), \(N_{\lambda}=30\), 3D}\label{fig:nbit_gmres_Robin_3D_vsnweak_plot}
    \end{subfigure}\caption{Number of iterations with respect to the number of
    subdomains \(\mJ\) for a domain size growing like \(\mJ^{1/d}\).  Disk (left) and sphere (right) of
    increasing radius. \textsc{GMRes}
    algorithm.}\label{fig:nbit_Robin_2D_vsnweak_plot}
\end{figure}

In the second set of results, the domain \(\Omega\) increases in size as the
number of subdomains \(\mJ\) grows so as to keep a fixed number of degrees of
freedom per subdomain.
Figure~\ref{fig:nbit_Robin_2D_vsnweak_plot} reports the iteration count with
respect to \(\mJ\) for the 2D and 3D cases.
The size of the domain is chosen to grow like \(\mJ^{1/d}\) where \(d\) is the
dimension of ambient space, so as to keep a fixed size (in terms of DOFs) for
the local subdomains.
In 2D the domain is a disk of radius increasing from \(R=1\) to \(R=16\), and
in 3D the domain is a sphere of radius increasing from \(R=1\) to \(R=3.7\).
The growth of the number of iteration to reach the set tolerance also appears
to scale like \(\mJ^{1/d}\).

\subsection{Heterogeneous medium}

\begin{figure}[h]
    \centering
    \includegraphics[width=0.6\textwidth]{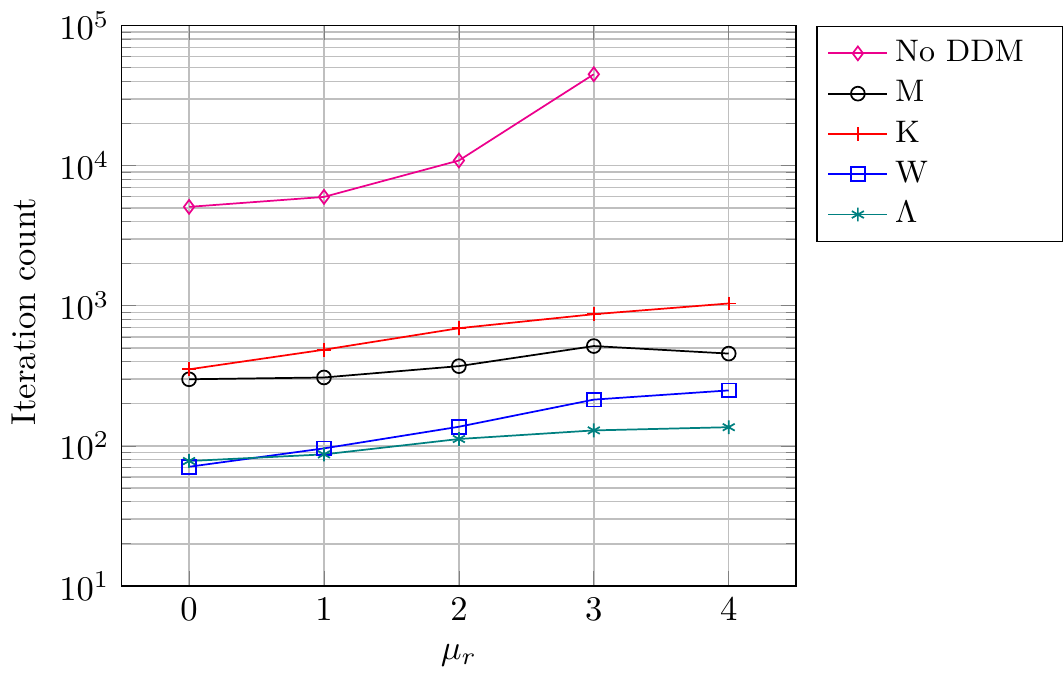}%
    \caption{Number of iterations with respect to increasing contrast in
    \(\mu\). Fixed parameters \(\kappa=10\), \(\mJ=10\), \(N_{\lambda}=50\),
    2D, disk of radius \(R=1\). \textsc{GMRes}
    algorithm.}\label{fig:nbit_Robin_2D_k10_n10_Nl50_vsmu_plot}
\end{figure}

\begin{figure}[h]
    \centering
    \includegraphics[width=0.6\textwidth]{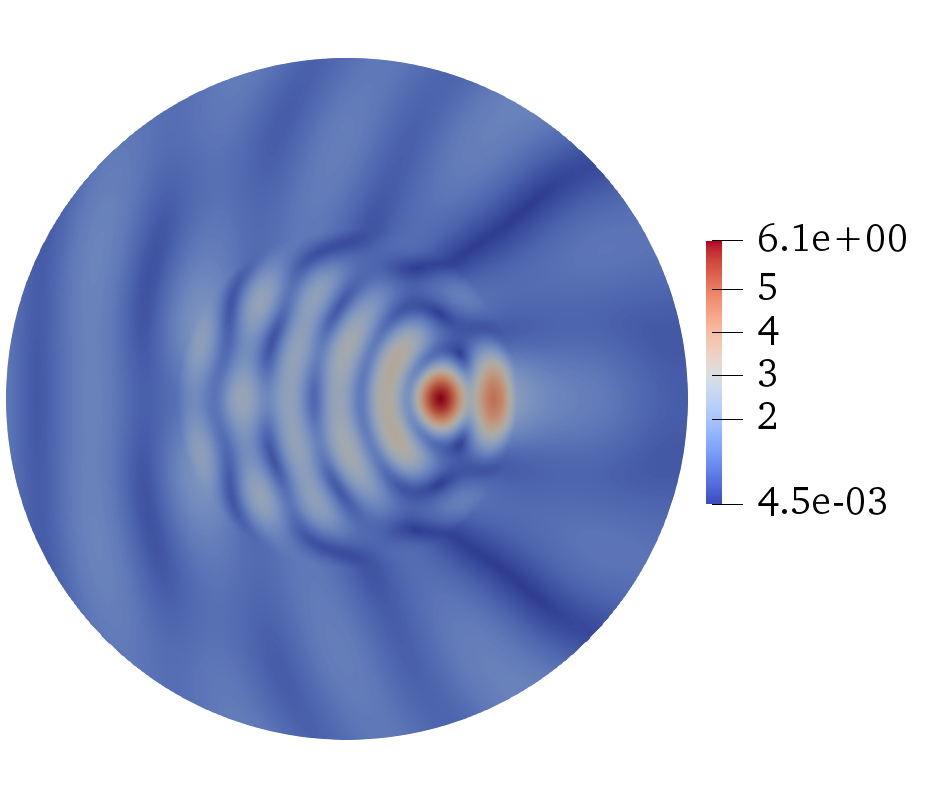}%
    \caption{Modulus of the total field. Fixed parameters \(\kappa=10\), \(N_{\lambda}=50\),
    2D, disk of radius \(R=1\), \(\mu_{r}=4\).}\label{fig:heterogenous_plot}
\end{figure}

We close the numerical experiment section with results in heterogeneous
medium in 2D.
The domain of propagation \(\Omega\) is still a disk of radius \(R=1\), but
this time with a circular inclusion of a different medium in the region with
radius \(R\leq 0.5\). The coefficient \(\mu\) is still equal to \(1\) outside
the inclusion and takes the value \(\mu=1+\mu_r\) inside, with \(\mu_{r}\)
varying from \(0\) (homogeneous medium) to \(4\).
Figure~\ref{fig:nbit_Robin_2D_k10_n10_Nl50_vsmu_plot} reports the iteration
counts for the \textsc{GMRes} algorithm as the medium varies.
The partition is composed of \(10\) subdomains so that some interfaces are cut
by the discontinuity in the medium.
One can observe that the number of iterations to get to convergence increases
greatly for the undecomposed problem (line plot labeled `No
DDM').
This is due to the appearance in the solution of "quasi-modes" of the
inclusion with large amplitude. For an illustration of this effect, the modulus
of the total field is represented in Figure~\ref{fig:heterogenous_plot} for the
value \(\mu_{r}=4\).
On the other hand, the DD algorithm performs well, with a
number of iterations only mildly growing.

\section*{Acknowledgments} This work was supported by the project NonlocalDD funded by
the French National Research Agency, grant ANR--15--CE23--0017--01. The authors would like to
thank Patrick Joly and Francis Collino for many inspiring discussions and, in
particular, for pointing a simplification in the proof of
Corollary~\ref{EstimCoercivite}. 

%\begin{acknowledgements}
%\end{acknowledgements}


\bibliographystyle{spmpsci}      % mathematics and physical sciences
\bibliography{biblio}   % name your BibTeX data base

%% Non-BibTeX users please use
%\begin{thebibliography}{}
%%
%% and use \bibitem to create references. Consult the Instructions
%% for authors for reference list style.
%\end{thebibliography}

\end{document}